\documentclass{article}
\usepackage{times,cite}
\usepackage[T1]{fontenc}
\usepackage[utf8]{inputenc}
\usepackage[a4paper, total={6.5in, 8in}]{geometry}
\usepackage{graphicx}
\usepackage{wrapfig}

\usepackage{amsfonts, amsthm, amsmath, amssymb, amstext, latexsym}
\usepackage{tikz}
\usepackage{pgfplots}
\pgfplotsset{
        table/search path={../../pics/}
}
\pgfplotsset{compat=1.17}
\usepackage{float} 

\usepackage{algorithm}
\usepackage{algpseudocode}

\usepackage{upgreek} 
\usepackage{comment}
\usepackage{cleveref}

\usepackage{caption}

\usepackage{siunitx}

\newtheorem{theorem}{Theorem}[section]

\newtheorem{remark}[theorem]{Remark}
\newtheorem{example}[theorem]{Example}

\usepackage{enumitem}
\newlist{todolist}{itemize}{2}
\setlist[todolist]{label=$\square$}
\usepackage{pifont}
\newcommand{\R}{\mathbb{R}}

\newcommand{\Lc}{\mathcal{L}}

\newcommand{\Jc}{\mathcal{J}}

\newcommand{\Ec}{\mathcal{E}}
\newcommand{\gradL}{\nabla\mathrm{L}(\uprho,\mathrm{u},\mathrm{p})}

\newcommand{\Norm}[1]{\|#1\|}
\newcommand{\diag}[1]{\text{diag(}#1\text{)}}

\newcommand{\diverg}[1]{\text{div}\left(#1\right)}

\newcommand{\dLn}[1]{\partial_{\mathrm{#1}} \mathrm{L}(\uprho,\mathrm{u},\mathrm{p})}
\newcommand{\dLninit}[1]{\partial_{\mathrm{#1}} \mathrm{L}(\uprho_0,\mathrm{u}_0,\mathrm{p}_0)}

\newcommand{\prox}[2]{\ifthenelse{\equal{#2}{}}{\text{prox}_{#1}}{\text{prox}_{#1}(#2)}}
\newcommand{\innerpr}[2]{\langle #1,#2 \rangle}

\newcommand{\ddLn}[2]{\partial_{\mathrm{#1}}\partial_{\mathrm{#2}} \mathrm{L}(\uprho,\mathrm{u},\mathrm{p})}

\title{Continuation methods for higher-order topology optimization}

\author{P. Gangl \\
  M. Winkler\\
  \emph{Johann Radon Institute for Computational and Applied }\\\emph{Mathematics (RICAM), Austrian Academy of Sciences}}

\begin{document}

\maketitle

\begin{abstract}
  We aim to solve a topology optimization problem where the distribution of material in the design domain is represented by a density function.
  To obtain candidates for local minima, we want to solve the first order optimality system via Newton's method. This requires the
  initial guess to be sufficiently close to the a priori unknown solution. Introducing a stepsize rule often allows for
  less restrictions on the initial guess while still preserving convergence. In topology optimization one typically encounters nonconvex problems
  where this approach might fail. We therefore opt for a homotopy (continuation) approach which is based on solving a sequence of parametrized problems to approach the solution of the original problem.
  In the density based framework the values of the design variable are constrained by $0$ from below and $1$ from above. Coupling the homotopy method with a barrier strategy enforces these constraints to be satisified.
  The numerical results for a PDE-constrained compliance minimization problem demonstrate that this combined approach maintains feasibility of the density function and
  converges to a (candidate for a) locally optimal design without a priori knowledge of the solution.
\end{abstract} 

\section{Introduction}
Topology optimization methods aim at finding an optimal distribution of different materials within a design region in order to minimize a certain quantity of interest, possibly while meeting some constraints. The most widely used classes of methods are density-based methods \cite{BendsoeSigmund2003} where the design is represented by a continuous density variable with values in the interval $[0,1]$ with a penalization of intermediate values; and level set methods where the sign of a continuous level set function determines the local material assignment \cite{allaire2020survey}. Topology optimization has been successfully applied to a wide variety of problems from different physical disciplines including mechanical engineering \cite{allaire2004structural, BendsoeSigmund1999} and electromagnetics \cite{Cherrire2022}, and even to coupled multiphysics problems \cite{GanglKrennDeGersem2025, FepponAllaireDapognyJolivet2020}.

The most widely used class of numerical optimization algorithms are based on first order derivative information and include (projected) gradient descent or sequential convex approximation methods such as the method of moving asymptotes (MMA) \cite{Svanberg1987}. While these methods have proven successful in practice, they may require a large number of iterations until convergence is reached. For that reason, also second order methods for topology optimization have been investigated, see e.g. \cite{AdamSurowiec2018} for an interior point method or \cite{Evgrafov2014} for a state space Newton method.
Beside their fast local convergence, Newton methods can be useful for applications where a solution path is to be followed using a predictor-corrector scheme. For instance, in \cite{PapadopoulosFarrellSurowiec2021}, multiple solution paths were encountered by a deflation approach and subsequently followed resulting in different local solutions to nonconvex topology optimization problems. As a second application that serves as motivation for our work, we mention the tracing of Pareto curves in multi-objective optimization \cite{schmidt2008pareto, MartinSchuetze2018}.

Since Newton-type methods are only locally convergent, globalization strategies are important.
In this work, we aim at finding a local solution to a topology optimization problem by determining a solution to the first order optimality system using Newton's method. In order to obtain convergence, we employ a homotopy method \cite{allgower2012numerical}. The idea of homotopy methods for solving systems of nonlinear equations $F(x)=0$ is to connect the problem of interest with a much simpler problem $G(x)=0$, whose solution is known or easy to obtain, by defining the family of problems
\begin{equation} \label{eq_defH}
    H(x,t)=t F(x)+(1-t)G(x).
\end{equation}
Then, starting from $t=0$ the homotopy parameter $t$ is gradually increased and the new problem \eqref{eq_defH} with the updated value of $t$ is solved by application of Newton's method starting out from the solution at the previous step. This approach has been successfully applied for finding solutions to challenging systems of nonlinear equations \cite{Malinen2010} or optimization problems \cite{Dunlavy2005} and has recently been applied in the context of shape optimization \cite{CesaranoEndtmayerGangl2024PRE}. In order to incorporate the box constraints on the density variable, we combine the homotopy approach with a primal-dual barrier method.

The remainder of this article is organized as follows: In Section \ref{sec:model_prb} we introduce the model problem from linear elasticity. An overview over homotopy and barrier methods is given in Sections \ref{sec:hmtpy} and \ref{sec:barrier} before combining the two concepts in Section \ref{sec:barrier_hmtpy}. Finally, the method is illustrated in numerical experiments in Section \ref{sec:results}.

\section{Model problem} \label{sec:model_prb}

We consider a rectangular computational domain $\Omega \subset \R^2$ with boundary $\partial \Omega=\Gamma_D \cup \Gamma_N$ and the subboundaries
$\Gamma_D = \Gamma_{D,0}$ and $\Gamma_N = \Gamma_{N,g_N} \cup \Gamma_{N,0}$.
The distribution of two materials on $\Omega$ is represented by the density variable $\rho:\Omega \to [0,1]$ with $\rho(x)=0$ and
$\rho(x)=1$ indicating absence and existence of a material at a given point, respectively.
We now apply a traction force $g_N$ on $\Gamma_{N,g_N}$, fix the boundary $\Gamma_D$ and aim to
distribute the material in such a way that the compliance of the mechanical structure is minimized or, equivalently, its stiffness maximized.
The resulting small displacement $u$ caused by the force $g_N$ based on a given material distribution $\rho$ is modelled by the linear elasticity equation which in its weak form is expressed by the state equation $e(\rho,u)=0$ with $\rho \in Q$
for a given design space $Q$ and $u \in V := (H_{0,\Gamma_D}^1(\Omega))^2$.
To exclude the trivial solution of full material we require the volume $\int_\Omega \rho \, dx$ to be constrained from above.
This is achieved by adding a volume control term to the objective weighted by a positive parameter $\gamma>0$.

We summarize the preceding discussion by stating the corresponding PDE constrained minimization problem in its reduced form,
\begin{align}
  \underset{\rho \in Q}{\mbox{inf}} \, j(\rho):=j(\rho,u(\rho)) \quad \text{subject to} \quad 0 \leq \rho \leq 1 \text{ a.e. on } \Omega
\end{align}
where $j(\rho,u):= \int_{\Gamma_{N,g_N}} g_N \cdot u \, ds + \gamma \int_\Omega \rho \, dx$ and $u=u(\rho) \in V$ solves the state equation $e(\rho,u)=0$ in $V^*$
which is given by
\begin{align} \label{eq_state}
  \innerpr{e(\rho,u)}{v}_{V^*\times V}:= \int_\Omega (\lambda^L(\rho) \diverg{u}Id + 2\mu^L(\rho)\Ec(u)) : \Ec(v) \, dx - \int_{\Gamma_{N,g_N}} g_N \cdot v \, ds
\end{align}
for $v \in V$. Here we denote by $\Ec(u):= \frac{1}{2}(\nabla u + \nabla u^\top)$ the linearized strain tensor and define the interpolation of the Lamé parameters based on the material distribution as $\lambda^L(\rho):=\lambda^L_0 + \rho^p(\lambda^L_1 - \lambda^L_0)$
and $\mu^L(\rho):= \mu^L_0 + \rho^p(\mu^L_1 - \mu^L_0)$ for given positive $\lambda^L_0,\lambda^L_1,\mu^L_0,\mu^L_1$ and $p \geq 1$ which we set to $p=3$ in the following.

In this formulation the optimization problem may not admit a solution which may result in numerical instabilities such as checkerboard patterns or mesh dependence of numerical solutions \cite{Sigmund1998}.
We choose the regularization approach based on the Ginzburg-Landau energy $P_\varepsilon^{0,1}(\rho):=\varepsilon\int_\Omega |\nabla \rho|^2 \, dx + \frac{1}{\varepsilon} \int_\Omega \rho(1-\rho) \, dx$ which can be shown to $\Gamma$-converge to the perimeter functional, $P^{0,1}_\varepsilon \overset{\Gamma}{\to} P^{0,1}$ as $\varepsilon \searrow 0$. More precisely, the limit $P^{0,1}(\rho)$ corresponds to the perimeter of ${\{\rho=1\}}$ associated with a binary $0$-$1$-restriction imposed on the values of $\rho$ \cite{Modica1987}.
Thus adding $P_\varepsilon^{0,1}(\rho)$ to the objective with a sufficiently large weight $\beta>0$ and a fixed positive $\varepsilon$ results in forcing the density variable to only take
values close to $0$ or $1$ (i.e., eliminating intermediate materials) and simultaneously controlling the overall perimeter avoiding checkerboard patterns.
Since this regularization involves the gradient of $\rho$, we require for the density space that $Q:=H^1(\Omega)$.

Eventually, the regularized infinite dimensional reduced PDE-constrained problem is given by
\begin{align}
    \underset{\rho \in Q}{\min} \, \Jc(\rho):=J(\rho,u(\rho)) \quad \text{subject to} \quad 0 \leq \rho \leq 1 \text{ a.e. on } \Omega
    \label{pr:red_abstract}
\end{align}
with $J(\rho,u):= j(\rho,u) + \frac{\beta}{2} P_\varepsilon^{0,1}(\rho)$ and $u(\rho)$ satisfying the state equation \eqref{eq_state} as before for a given design $\rho$. A discretization of state and design variable by means of piecewise linear and globally continuous finite elements results in a discretized version of \eqref{pr:red_abstract} which is of the form
\begin{align} \label{eq_abstract_discr}
  \underset{x \in \R^n}{\min} \, f(x) \quad \text{subject to} \quad a \leq x \leq b
\end{align}
where $a,b \in \R^n$ such that $a<b$ componentwise. Here, the vector $x$ represents the finite element coefficient vector $\uprho$ with respect to the hat basis functions of $\rho_h(y) = \sum_{i=1}^n \uprho_i \varphi_i(y) \approx \rho(y)$ for $y \in \Omega$.
The subsequent discussion is concerned with computing stationary points of the abstract discretized problem \eqref{eq_abstract_discr}.

\section{Homotopy method}{\label{sec:hmtpy}}

In this section we temporarily ignore the constraints $a \leq x \leq b$ and are concerned instead
with computing stationary points of the \textit{unconstrained} minimization problem $\underset{x \in \R^n}{\min} \, f(x)$
with differentiable but nonconvex objective $f:\R^n \to \R$.
This may be achieved by finding a point $x^*$ satisfying the first order necessary optimality condition, i.e., $\nabla f(x^*)=0$.

More abstractly, we consider a system of equations $F:\R^n \to \R^n, \, F(x)=0$ with $F$ being continuously differentiable but nonconvex.
Computing a root via Newton's method requires the initial guess to already be close to a solution $x^*$.
This a priori knowledge cannot be expected in general. Coupling Newton's method with a linesearch strategy as a way
of enlargening the basin of attraction may also fail due to certain assumptions on the regularity of the Jacobian of $F$ that may be violated in the nonconvex case \cite{Schwetlick1979}.

The homotopy method avoids the requirements of proximity via modifiying the target
problem $F(x)=0$ by adding a perturbation term $G(x)$ with $G:\R^n \to \R^n$ continuously differentiable whose root is known or can easily be computed and is further
defined as the initial guess $x_0$. In the course of the method the contribution of the perturbation term is
gradually decreased until the target problem is reached.
In other words, a smooth transition (continuation) between the auxiliary problem $G(x)=0$ and the target problem $F(x)=0$
is required to be provided by the homotopy map $H:\R^n \times \R \to \R^n, (x,t) \mapsto H(x,t)$. The scalar homotopy parameter $t$ acts
as weight between the two problems such that $H(x,0)=G(x)$ and $H(x,1)=F(x)$.

A popular choice for constructing the homotopy map is the \textit{convex} homotopy
\begin{align}
  H(x,t):= t F(x) + (1-t) G(x)
\end{align}
where the mapping $G:\R^n \to \R^n$ may be chosen as $G(x):=x-x_0$ or $G(x):=F(x)-F(x_0)$.
In the latter case the resulting homotopy map is called the \textit{global homotopy} and reads as
\begin{align}
  H(x,t) = t F(x) + (1-t) (F(x)-F(x_0)) = F(x) - (1-t) F(x_0).
  \label{def:global_hmtpy}
\end{align}
We restrict ourselves to this particular choice in the following.

The following example illustrates the idea of the approach.
\begin{example}
  We consider the problem of computing a root of the equation $F(x) := 4x^3- 3x^2 - 2x + 1 = 0$. The choice $x_0 = -1.2$ induces
   the global homotopy map $H(x,t)=F(x)-(1-t)F(x_0)= F(x) - (1-t)c$ with $c:=-7.8319$. The perturbation of the target problem
  thus causes a translation parametrized by $t$.
  As $t$ approaches $1$ the translation constant goes to zero causing $H(x,t)$ to morph into the target function.
  The zeros of the homotopy map which move towards the root $x^*$ of $F$ are depicted in \Cref{fig:hmtpy_ex} for different parameters $t$.

\begin{figure}[H]
\begin{minipage}{85mm}
  \begin{figure}[H]
	\centering
	\scalebox{0.9}{%
	\begin{tikzpicture}[
          declare function={
            f(\x) = 4*\x*\x*\x - 3*\x*\x - 2*\x + 1;
            h(\x,\init,\t) = f(\x) - (1-\t)*f(\init);
            }
        ]
		\begin{axis}[axis lines=middle,
            axis on top = false,
            xmin = -2, xmax = 2,
			ymax=15,ymin=-15,xtick=\empty,ytick=\empty,
            xlabel={$x$},
            ticklabel style={above left},
			smooth,samples=201,clip=true]

            \addplot[red, domain=-2:2, smooth,thick] {h(x,-1.2,0)} node[left,pos=0.79] {$H(x,0)$};
            \addplot[blue, domain=-2:2, smooth,thick] {f(x)} node[right,pos=0.5] {$F(x)$};

            \filldraw [gray] (-1.2,{h(-1.2,-1.2,0)}) circle (1.5pt) node[above left] {$x(0)$};
            \filldraw [black] (-0.6404,0) circle (1.5pt) node[below right] {$x^*$};

            \draw[<->,line width = 0.3pt] (0.2,{f(0.2)}) -- (0.2,{h(0.2,-1.2,0)}) node[right,pos=0.5] {$c$};
			\end{axis}
		\end{tikzpicture}}
	\end{figure}
\end{minipage}
\hfil
\begin{minipage}{85mm}
  \begin{figure}[H]
	\centering
	\scalebox{0.9}{%
	\begin{tikzpicture}[
          declare function={
            f(\x) = 4*\x*\x*\x - 3*\x*\x - 2*\x + 1;
            h(\x,\init,\t) = f(\x) - (1-\t)*f(\init);
            }
        ]
		\begin{axis}[axis lines=middle,
            axis on top = false,
            xmin = -2, xmax = 2,
			ymax=15,ymin=-15,xtick=\empty,ytick=\empty,
            xlabel={$x$},
            ticklabel style={above left},
			smooth,samples=201,clip=true]

            \addplot[red, domain=-2:2, smooth,thick] {h(x,-1.2,0.4)} node[left,pos=0.79] {$H(x,0.4)$};
            \addplot[blue, domain=-2:2, smooth,thick] {f(x)} node[right,pos=0.5] {$F(x)$};

            \filldraw [gray] (-1.0420,0) circle (1.5pt) node[above left] {$x(0.4)$};
            \filldraw [black] (-0.6404,0) circle (1.5pt) node[below right] {$x^*$};
			\end{axis}
		\end{tikzpicture}}
	\end{figure}
\end{minipage}
\end{figure}

\begin{figure}[H]
\begin{minipage}{85mm}
  \begin{figure}[H]
	\centering
	\scalebox{0.9}{%
	\begin{tikzpicture}[
          declare function={
            f(\x) = 4*\x*\x*\x - 3*\x*\x - 2*\x + 1;
            h(\x,\init,\t) = f(\x) - (1-\t)*f(\init);
            }
        ]
		\begin{axis}[axis lines=middle,
            axis on top = false,
            xmin = -2, xmax = 2,
			ymax=15,ymin=-15,xtick=\empty,ytick=\empty,
            xlabel={$x$},
            ticklabel style={above left},
			smooth,samples=201,clip=true]

            \addplot[red, domain=-2:2, smooth,thick] {h(x,-1.2,0.65)} node[left,pos=0.79] {$H(x,0.65)$};
            \addplot[blue, domain=-2:2, smooth,thick] {f(x)} node[right,pos=0.5] {$F(x)$};

            \filldraw [gray] (-0.9147,0) circle (1.5pt) node[above left] {$x(0.65)$};
            \filldraw [black] (-0.6404,0) circle (1.5pt) node[below right] {$x^*$};
			\end{axis}
		\end{tikzpicture}}
	\end{figure}
\end{minipage}
\hfil
\begin{minipage}{85mm}
  \begin{figure}[H]
	\centering
	\scalebox{0.9}{%
	\begin{tikzpicture}[
          declare function={
            f(\x) = 4*\x*\x*\x - 3*\x*\x - 2*\x + 1;
            h(\x,\init,\t) = f(\x) - (1-\t)*f(\init);
            }
        ]
		\begin{axis}[axis lines=middle,
            axis on top = false,
            xmin = -2, xmax = 2,
			ymax=15,ymin=-15,xtick=\empty,ytick=\empty,
            xlabel={$x$},
            ticklabel style={above left},
			smooth,samples=201,clip=true]

            \addplot[red, domain=-2:2, smooth,thick] {h(x,-1.2,0.9)} node[left,pos=0.79] {$H(x,0.9)$};
            \addplot[blue, domain=-2:2, smooth,thick] {f(x)} node[right,pos=0.5] {$F(x)$};

            \filldraw [gray] (-0.7399,0) circle (1.5pt) node[above left] {$x(0.9)$};
            \filldraw [black] (-0.6404,0) circle (1.5pt) node[below right] {$x^*$};
			\end{axis}
		\end{tikzpicture}}
	\end{figure}
\end{minipage}
\caption{The homotopy map for $t=0, 0.4, 0.65$ and $0.9$ in contrast to the graph of $F$}
\label{fig:hmtpy_ex}
\end{figure}
\end{example}

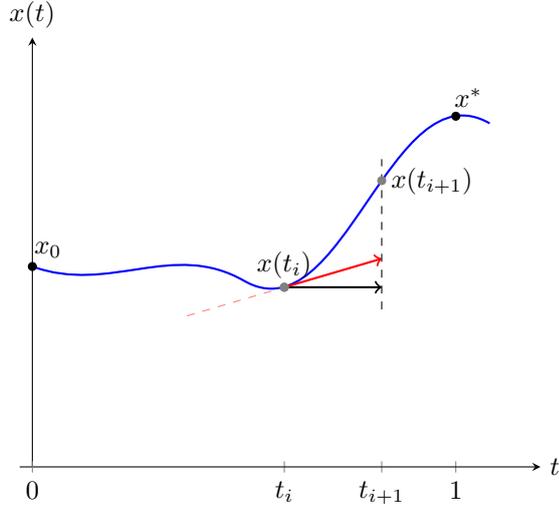
\begin{figure}
\begin{minipage}{75mm}
  \begin{figure}[H]
        \scalebox{1}{%
        \begin{tikzpicture}
            \begin{axis}[axis on top = false,
                axis y line=middle,
                axis x line=middle,
                xmin=-0.3,xmax=12,ymax=6,ymin=0,
                xtick=\empty,ytick=\empty,
                xlabel={$t$},
                ylabel={$x(t)$},
                x label style={anchor=west},
                y label style={anchor=south},
                xtick={5.95,5.95 + 2*1.15,10},
                xticklabels={$t_i$,$t_{i+1}$,$1$},
                extra x ticks = {0},
                extra x tick labels = {$0$}]

                \coordinate (X) at (10,4.9);
                \coordinate (Xn) at (10.3,4.9);

                \coordinate (A1) at (0,2.8);
                \coordinate (A2) at (5,2.6);
                \coordinate (A3) at (10.8,4.8);

                \coordinate (P1) at (5.95,2.51);
                \coordinate (P2) at (5.95 + 2*1.15,2.51 + 2*0.2);
                \coordinate (P0) at (5.95 - 2*1.15,2.51 - 2*0.2);

                \coordinate (Xt) at (P2 |- 4,4);

                \coordinate (X0) at (0,2.8);
                \coordinate (X0n) at (0.37,2.8);

                \draw [thick,color=blue] (A1) to[out=-20,in=150] (A2) to[out=-30,in=150] (A3);

                \draw[thick,->] (P1) -- (Xt |- P1);

                \draw[thick,->,red] (P1) -- (P2);
                \draw[dashed,-,opacity=0.5,red] (P0) -- (P1);
                \draw[dashed,-] (P2 |- 7.1,2.2) -- (P2 |- 7.1,4.3);

                \filldraw [gray] (P1) circle (1.5pt);

                \filldraw [black] (X) circle (1.5pt);
                \node at (Xn) [above] {$x^*$};

                \filldraw [black] (X0) circle (1.5pt);
                \node at (X0n) [above] {$x_0$};

                \filldraw [gray] (P1) circle (1.5pt);
                \node at (P1) [above] {$x(t_i)$};

                \filldraw [gray] (Xt) circle (1.5pt);
                \node at (Xt) [right] {$x(t_{i+1})$};
			\end{axis}
      \end{tikzpicture}}










  \caption{Visualization of tracing the zero curve of $H$}
  \label{fig:trace_h_curve}
  \end{figure}
\end{minipage}
\hspace{1.6cm}
\begin{minipage}{66mm}
  \begin{center}
  \begin{algorithm}[H]
  \caption{Homotopy method}\label{alg:pc_hmtpy}
    \begin{algorithmic}
      \Require $x_0 \in \R^n$, update rule $\theta:[0,1] \to [0,1]$
      \Ensure $x^* \in \R^n$ with $F(x^*)=0$
      \State $x \gets x_0$, $t \gets 0$
      \While{$t < 1$}
        \State $t \gets \theta(t)$
        \While{$\Norm{H(x,t)}>tol$}
          \State Solve $H_x(x,t) \Delta x = -H(x,t)$
          \State $x \gets x + \Delta x$
        \EndWhile
      \EndWhile
    \end{algorithmic}
  \end{algorithm}
\end{center}
\end{minipage}
\end{figure}

The path consisting of all tuples $(x,t)$ with $x=x(t)$ and $t$ satisfying $H(x,t)=0$ is denoted the zero curve of the homotopy map $H$.
Hence to compute a solution of $F(x)=0$ by means of the homotopy method the zero curve of
the homotopy equation $H(x,t)= 0$ is numerically traced starting from $(x(0),0)=(x_0,0)$ until $(x(1),1)=(x^*,1)$ is reached.
More precisely, from a given point $(x(t_i),t_i)$ on the curve the next point may be obtained by updating $t_i$ to $t_{i+1}=\theta(t_i)$ via an appropriate
update rule $\theta:[0,1] \to [0,1]$, e.g. $\theta(t_i)=\text{min}(t_i + \Delta t,1)$ for a fixed or adaptively chosen increment $\Delta t$.
Subsequently, the corrector equation $H(x,t_{i+1})=0$ for fixed $t_{i+1}$ is solved starting out from the initial guess given by the
(zero order) predictor $(x(t_i),t_{i+1})$ (see \Cref{fig:trace_h_curve} with the predictor represented by the black arrow).

The main assumption for the existence of a smooth zero curve is that $0$ is a regular value of $H(x,t)$ for
all $t \in \R$, i.e., for any given $t$ the Jacobian $H_x(x,t)$ is regular at all points $x$ satisfying $H(x,t)=0$ \cite{Schwetlick1979,allgower2012numerical}.
An immeditate consequence is that the solution curve may be parametrized with respect to the homotopy parameter $t$.
Further, by a continuity argument Newton's method is well defined for points near the zero curve and thus can be applied to solve the
corrector equation.The so-obtained predictor-corrector scheme embedded in the homotopy method is summarized in \Cref{alg:pc_hmtpy}.

\begin{remark} {\label[remark]{rem:h_predictor}}
  For a sufficiently small stepsize $\Delta t = \theta(t)-t$ the zero order predictor $(x,t + \Delta t)$ is close to the homotopy curve causing Newton's method for
  solving the corrector equation to converge.
  To enhance the accuracy of the predictor (and thus allowing for a larger stepsize), the tuple $(x,t + \Delta t)$ may be replaced by the first order (or tangent)
  predictor $(\tilde{x},t + \Delta t)$
  where $\tilde{x}=x+ \Delta t \, x'(t)$ with $x'(t) \in \R^n$ satisfying the equation
  \begin{align} \label{eq:h_predictor}
     H_x(x,t) x'(t) = - H_t(x,t).
  \end{align}
  The predictor equation \eqref{eq:h_predictor} is derived from differentiating $H(x(t),t)=0$ with respect to $t$.
  As depicted in \Cref{fig:trace_h_curve}, the first order predictor (represented by the red arrow) demonstrates a better approximation of the next point on the zero curve compared to the zero order predictor.
  In general, a good balance between a small step size (security of convergence) and a large step size (fast progress on the curve) is
  crucial for the design of an efficient algorithm.

\end{remark}

\section{Barrier method}\label{sec:barrier}

Reintroducing the constraints on $x$ we now aim to apply the homotopy approach to compute stationary points of \eqref{eq_abstract_discr}.
The barrier approach (or interior point approach) allows for reformulating the constrained problem \eqref{eq_abstract_discr} as
a parametrized sequence of unconstrained problems \cite{Forsgren2002IMN}, thus making the respective first order optimality conditions form a system of equations and therefore
ensuring compatibility with the homotopy approach. We proceed by introducing the barrier concept for an implicit handling of inequality constraints,
stating the continuation method (with respect to the associated barrier parameter) induced by it and eventually improving the efficiency of
the latter by incorporating the dual variable.

Firstly we note that problem \eqref{eq_abstract_discr} can be expressed as
\begin{align}
  \underset{x \in \R^n}{\min} \, f(x) \quad \text{subject to}
  \quad
  c(x) \geq 0
  \label{pr:n_dim_box_expl}
\end{align}
where $c:\R^n \to \R^m$. The case of box constraints $a \leq x \leq b$ can be treated by $c(x):=\begin{bmatrix}
  c^a(x)^\top,c^b(x)^\top
\end{bmatrix}^\top$ where $m=2n$ and $c^a(x):=x-a$ and $c^b(x):=b-x$. We denote by $\mathcal{F}:=\{x \in \R^n: \, c(x)\geq 0\}$
the feasible region of the minimization problem \eqref{pr:n_dim_box_expl}.

The explicit constraints are further replaced by adding a barrier term $I(x)$ with $I:\R^n \to \R$ to the
objective where the contribution of said term to the objective is controlled by a positive weight, i.e., the barrier parameter $\mu > 0$.
We require the barrier term to be differentiable in the interior of the feasible domain $\mathcal{F}$ and that $I(x) \to + \infty$ for any strictly feasible
$x$ approaching a point on the boundary of $\mathcal{F}$.
Consequently, for sufficiently large $\mu$ the (unconstrained) minimizer of the barrier function $B(x;\mu):=f(x) + \mu I(x)$ is located
in the interior of the feasible region $\mathcal{F}$.
Here we opt for the logarithmic barrier term $I(x):= - \sum_{i=1}^m \, \log(c_i(x))$.

Implicitly assuming $c_i(x)>0$ for $i=1,...,m$ we may now reformulate \eqref{pr:n_dim_box_expl} for a given barrier parameter $\mu >0$ as
\begin{align}
  \underset{x \in \R^n}{\min} \, B(x;\mu):=f(x) - \mu \sum_{i=1}^{m} \log(c_i(x))
\end{align}
and state the corresponding necessary first order optimality condition
\begin{align}
  \nabla_x B(x;\mu)= \nabla f(x) - \mu \sum_{i=1}^{m} \frac{1}{c_i(x)} \nabla c_i(x) = 0.
\end{align}
These equations are then solved for $x$ with Newton's method for a decreasing sequence of positive parameters $(\mu_k)_k$ starting with $\mu_0 \gg 0$ with each solution $x(\mu_k)$ acting as an
initial guess for the updated problem $\nabla_x B(x;\mu_{k+1})=0$ until the target value $\mu_\infty \approx 0$ is reached.

We illustrate the procedure in a simple example.
\begin{example}
  The function $f(x) := x^4 - x^3 - x^2 + x + 0.25$ possesses two local minimizers within the interval $[a,b]=[-0.5,1]$.
  To compute a stationary point while remaining in $[a,b]$ the barrier function $B(x;\mu) := x^4 - x^3 - x^2 + x + 0.25 - \mu(\log(x+0.5)+\log(1-x))$
  is minimized successively. \Cref{fig:barrier_visual} depicts that the minimum of the barrier function approaches a constrained minimum of $f$ as $\mu$ is decreased.

\begin{figure}[H]
\begin{minipage}{80mm}
  \begin{figure}[H]
				\centering
				\scalebox{1}{%
				\begin{tikzpicture}[
          declare function={
            f(\x,\mu) = \x*\x*\x*\x - \x*\x*\x - \x*\x + \x + 0.25 - \mu*(ln(\x+0.5)+ln(1-\x));
            }
        ]
					\begin{axis}[axis lines=middle,
            axis on top = false,
						domain=-10:10,ymax=3.5,ymin=-0.5,xtick=\empty,ytick=\empty,
            xtick={-0.5,1},
            xticklabels={$a$,$b$},
            ticklabel style={above left},
						smooth,samples=201,clip=true]
            \addplot[blue, domain=-1.2:-0.5, smooth,thick,dashed] {f(x,0)};
						\addplot[blue, domain=-0.5:1, smooth,thick] {f(x,0)};
            \addplot[blue, domain=1:1.5, smooth,thick,dashed] {f(x,0)} node[right,pos=0.2] {$f(x)$};
            \addplot[red, domain=-0.5+0.01:1-0.01, smooth,thick] {f(x,2.9)} node[right,pos=0.5] {$B(x;2.9)$};
            \draw[dashed,-] (-0.5,-0.5) -- (-0.5,3.5);
            \draw[dashed,-] (1,-0.5) -- (1,3.5);

            \filldraw [gray] (0.2008,{f(0.2008,2.9)}) circle (1.5pt);
            \filldraw [black] (-0.5,{f(-0.5,0)}) circle (1.5pt) node[right,,xshift=0.1cm] {$f(x^*)$};
				 	\end{axis}
				\end{tikzpicture}}
	\end{figure}
\end{minipage}
\hfil
\begin{minipage}{80mm}
  \begin{figure}[H]
				\centering
				\scalebox{1}{%
				\begin{tikzpicture}[
          declare function={
            f(\x,\mu) = \x*\x*\x*\x - \x*\x*\x - \x*\x + \x + 0.25 - \mu*(ln(\x+0.5)+ln(1-\x));
            }
        ]
					\begin{axis}[axis lines=middle,
            axis on top = false,
						domain=-10:10,ymax=3.5,ymin=-0.5,xtick=\empty,ytick=\empty,
            xtick={-0.5,1},
            xticklabels={$a$,$b$},
            ticklabel style={above left},
						smooth,samples=201,clip=true]
            \addplot[blue, domain=-1.2:-0.5, smooth,thick,dashed] {x^4 - x^3 - x^2 + x + 0.25};
						\addplot[blue, domain=-0.5:1, smooth,thick] {x^4 - x^3 - x^2 + x + 0.25};
            \addplot[blue, domain=1:1.5, smooth,thick,dashed] {x^4 - x^3 - x^2 + x + 0.25} node[right,pos=0.2] {$f(x)$};
            \addplot[red, domain=-0.5+0.01:1-0.01, smooth,thick] {x^4 - x^3 - x^2 + x + 0.25 - 1.1*(ln(x+0.5)+ln(1-x))} node[left,pos=0.6] {$B(x;1.1)$};
            \draw[dashed,-] (-0.5,-0.5) -- (-0.5,3.5);
            \draw[dashed,-] (1,-0.5) -- (1,3.5);

            \filldraw [gray] (0.0315,{f(0.0315,1.1)}) circle (1.5pt);
            \filldraw [black] (-0.5,{f(-0.5,0)}) circle (1.5pt) node[right,xshift=0.1cm] {$f(x^*)$};
				 	\end{axis}
				\end{tikzpicture}}
	\end{figure}
\end{minipage}
\end{figure}

\begin{figure}[H]
\begin{minipage}{80mm}
  \begin{figure}[H]
				\centering
				\scalebox{1}{%
				\begin{tikzpicture}[
          declare function={
            f(\x,\mu) = \x*\x*\x*\x - \x*\x*\x - \x*\x + \x + 0.25 - \mu*(ln(\x+0.5)+ln(1-\x));
            }
        ]
					\begin{axis}[axis lines=middle,
            axis on top = false,
						domain=-10:10,ymax=3.5,ymin=-0.5,xtick=\empty,ytick=\empty,
            xtick={-0.5,1},
            xticklabels={$a$,$b$},
            ticklabel style={above left},
						smooth,samples=201,clip=true]
            \addplot[blue, domain=-1.2:-0.5, smooth,thick,dashed] {x^4 - x^3 - x^2 + x + 0.25};
						\addplot[blue, domain=-0.5:1, smooth,thick] {x^4 - x^3 - x^2 + x + 0.25};
            \addplot[blue, domain=1:1.5, smooth,thick,dashed] {x^4 - x^3 - x^2 + x + 0.25} node[right,pos=0.2] {$f(x)$};
            \addplot[red, domain=-0.5+0.01:1-0.01, smooth,thick] {x^4 - x^3 - x^2 + x + 0.25 - 0.4*(ln(x+0.5)+ln(1-x))} node[above,pos=0.5] {$B(x;0.4)$};
            \draw[dashed,-] (-0.5,-0.5) -- (-0.5,3.5);
            \draw[dashed,-] (1,-0.5) -- (1,3.5);

            \filldraw [gray] (-0.2456,{f(-0.2456,0.4)}) circle (1.5pt);
            \filldraw [black] (-0.5,{f(-0.5,0)}) circle (1.5pt) node[right,xshift=0.1cm] {$f(x^*)$};
				 	\end{axis}
				\end{tikzpicture}}
	\end{figure}
\end{minipage}
\hfil
\begin{minipage}{80mm}
  \begin{figure}[H]
				\centering
				\scalebox{1}{%
				\begin{tikzpicture}[
          declare function={
            f(\x,\mu) = \x*\x*\x*\x - \x*\x*\x - \x*\x + \x + 0.25 - \mu*(ln(\x+0.5)+ln(1-\x));
            }
        ]
					\begin{axis}[axis lines=middle,
            axis on top = false,
						domain=-10:10,ymax=3.5,ymin=-0.5,xtick=\empty,ytick=\empty,
            xtick={-0.5,1},
            xticklabels={$a$,$b$},
            ticklabel style={above left},
						smooth,samples=201,clip=true]
            \addplot[blue, domain=-1.2:-0.5, smooth,thick,dashed] {x^4 - x^3 - x^2 + x + 0.25};
						\addplot[blue, domain=-0.5:1, smooth,thick] {x^4 - x^3 - x^2 + x + 0.25};
            \addplot[blue, domain=1:1.5, smooth,thick,dashed] {x^4 - x^3 - x^2 + x + 0.25} node[right,pos=0.2] {$f(x)$};
            \addplot[red, domain=-0.5+0.01:1-0.01, smooth,thick] {x^4 - x^3 - x^2 + x + 0.25 - 0.1*(ln(x+0.5)+ln(1-x))} node[above,pos=0.72] {$B(x;0.1)$};
            \draw[dashed,-] (-0.5,-0.5) -- (-0.5,3.5);
            \draw[dashed,-] (1,-0.5) -- (1,3.5);

            \filldraw [gray] (-0.41,{f(-0.41,0.1)}) circle (1.5pt);
            \filldraw [black] (-0.5,{f(-0.5,0)}) circle (1.5pt) node[right,xshift=0.1cm] {$f(x^*)$};
				 	\end{axis}
				\end{tikzpicture}}
	\end{figure}
\end{minipage}
\caption{The barrier function for different parameters $\mu=2.9,1.1,0.4,0.1$}
\label{fig:barrier_visual}
\end{figure}
\end{example}

Two main difficulties arise when emplyoing this approach.
Firstly, the Hessian $\nabla_x^2 B(x;\mu)$ is ill-conditioned for small $\mu$ and $x$ close to the solution of \eqref{pr:n_dim_box_expl}.
This however does not affect the accuracy of the computed Newton direction significantly \cite{Forsgren2002IMN}.
Secondly, the radius of the basin of attraction for Newton's method tends to zero as $\mu$ approaches zero, thus allowing only very small steps
of the barrier parameter to be taken and therefore making this method very inefficient \cite{Villalobos2004}.

At the cost of introducing the dual variable $z:= \mu \oslash c(x) \in \R^m$ with $\oslash$ denoting the componentwise division, i.e., $(\mu \oslash c(x))_i = \mu / c_i(x)$, $i=1, \dots, m$, (and
thus adding $m$ additional additional equations to the optimality system) these deficiencies can be alleviated leading to the primal-dual barrier formulation.
The definition of the dual variable is directly incorporated in the optimality system resulting in
\begin{align}
  F(x,z;\mu) := \begin{bmatrix}
    \nabla f(x) - (Jc(x))^\top z \\
    z \odot c(x) - \mu e_m
  \end{bmatrix}
  =
  \begin{bmatrix}
    0_n \\
    0_m
  \end{bmatrix}
  \label{eq:pd_system}
\end{align}
  where
\begin{align}
  Jc(x) := \begin{bmatrix}
    - \nabla c_1(x)^\top -\\
    \vdots \\
    - \nabla c_m(x)^\top -
  \end{bmatrix} \in \R^{m \times n},
  \,
  e_m =
  \begin{bmatrix}
    1 \\
    \vdots \\
    1
  \end{bmatrix}
  \in \R^m
\end{align}
and $\odot$ denoting the componentwise multiplication of two column vectors, i.e., for $a, b \in \mathbb R^m$, the componentwise product $a \odot b \in \mathbb R^m$ is defined by $(a \odot b)_i = a_i b_i$, $i=1, \dots, m$.

\begin{remark}
  The system \eqref{eq:pd_system} may also be derived by perturbing the complementarity conditions of the KKT-system of an inequality constrained minimization problem.
  Thus the dual variable $z$ tends to the associated Lagrange multiplier $\lambda$ as the barrier parameter $\mu$ approaches zero \cite{Forsgren2002IMN}.
\end{remark}

For a given $\mu >0$, the problem \eqref{eq:pd_system} may be solved for the tuple $(x,z) \in \R^n \times \R^m$ by applying Newton's method with each update step $(\Delta x, \Delta z) \in \R^n \times \R^m$
satisfying
\begin{align}
  \begin{bmatrix}
    \nabla^2 f(x) - \sum_{i=1}^m z_i \nabla^2 c_i(x) & - (Jc(x))^\top \\
    Z Jc(x) & C(x)
  \end{bmatrix}
  \begin{bmatrix}
    \Delta x \\
    \Delta z
  \end{bmatrix}
  =
  -
  \begin{bmatrix}
    \nabla f(x) - (Jc(x))^\top z \\
    z \odot c(x) - \mu e_m
  \end{bmatrix}
  \label{eq:pd_newton}
\end{align}

where $Z:= \diag{z_i}$ and $C(x):= \diag{c_i(x)}$.
The corresponding radius of the sphere of convergence is indeed bounded from below by a strictly positive constant as $\mu$ tends to zero \cite{Villalobos2004}.

In its fundamental form the primal-dual barrier method is stated in \Cref{alg:pd_barrier}.
Note that the algorithm only gives an approximation of the solution of the KKT-system
of \eqref{pr:n_dim_box_expl} with its quality increasing the closer to zero the target parameter $\mu_\infty$ is chosen.

\begin{figure}
  \centering
  \begin{minipage}{0.6\linewidth}
    \begin{algorithm}[H]
    \caption{Primal dual barrier method}\label{alg:pd_barrier}
      \begin{algorithmic}
        \Require $x_0 \in \R^n$ strictly feasible, update rule $\theta: \R \to \R$, $\mu_0 \gg 0$, $\mu_\infty \approx 0$
        \Ensure $(x,z)$ approximately fulfilling the KKT-system of \eqref{pr:n_dim_box_expl}
        \State $x \gets x_0$, $\mu \gets \mu_0$
        \State $z \gets \mu_0 \oslash c(x_0)$
        \While{$\mu \geq \mu_\infty$}
          \State $\mu \gets \theta(\mu)$
          \While{$\Norm{F(x,z;\mu)}>tol$}
            \State Solve $JF(x,z;\mu) \begin{bmatrix}
              \Delta x \\
              \Delta z
            \end{bmatrix}=
            -F(x,z;\mu)$ as in \eqref{eq:pd_newton}
            \State $x \gets x + \Delta x$
            \State $z \gets z + \Delta z$
          \EndWhile
        \EndWhile
      \end{algorithmic}
    \end{algorithm}
  \end{minipage}
\end{figure}


\begin{remark}
  The initial guess $x_0$ is preferably chosen in the analytic center of the feasible region, i.e., as the minimizer of $(-1) \cdot\sum_{i=1}^{m} \log(c_i(x))$.
  For updating the barrier parameter a simple contraction approach may be pursued, i.e., $\theta(\mu)=\alpha \cdot \mu$ for some $\alpha \in (0,1)$. More advanced techniques
  addressing this aspect can be found in \cite{Nocedal2006NOpt}.
\end{remark}

 In the case of only box constraints setting $z=[(z^a)^\top,(z^b)^\top]^\top \in \R^{2n}$ in \eqref{eq:pd_system} and \eqref{eq:pd_newton} yields
 \begin{align}
    F_{box}(x,z^a,z^b, \mu) :=
    \begin{bmatrix}
        \nabla f(x) - z^a + z^b \\
        z^a \odot c^a(x) - \mu e_n \\
        z^b \odot c^b(x) - \mu e_n
    \end{bmatrix}
    =
    \begin{bmatrix}
      0_n \\
      0_n \\
      0_n
    \end{bmatrix}
  \end{align}
  and
  \begin{align}
    \begin{bmatrix}
    \nabla^2 f(x) & -I_n &I_n \\
    Z^a & C^a(x) & 0 \\
    -Z^b & 0 & C^b(x)
  \end{bmatrix}
  \begin{bmatrix}
    \Delta x \\
    \Delta z^a \\
    \Delta z ^b
  \end{bmatrix}
  =
  -
  \begin{bmatrix}
    \nabla f(x) - z^a + z^b \\
     z^a \odot c^a(x) - \mu e_n \\
    z^b \odot c^b(x) - \mu e_n
  \end{bmatrix}
  \label{eqn:pd_newton}
 \end{align}
where $Z^a:= \diag{z^a_i}, Z^b:= \diag{z^b_i}$ and $C^a(x):= \diag{c^a_i(x)}$, $C^b(x):= \diag{c^a_i(x)}$ accordingly.
The analytic center of the feasible region (and therefore the initial guess) can be explicitly computed in this setting
and is given by $x_0:=\frac{a+b}{2}$.

\section{A barrier homotopy approach}{\label{sec:barrier_hmtpy}}

We derive a combined approach exploiting both the robust convergence properties of the homotopy method (in terms of not requiring proximity of the initial guess to the solution)
and the implicit handling of the inequality constraints (here: box constraints) by the barrier method.

To this end, we proceed by recalling the global homotopy map $H(x,t):= F(x) - (1-t) F(x_0)$.
We further suppose that $F:\R^n \to \R^n$ represents an optimality function for the objective $f:\R^n \to \R$,
i.e., $F(x)=0$ defines a necessary condition for $x$ minimizing $f$.
In the unconstrained case the function $F$ coincides with the gradient of $f$.
For a minimization problem with $m$ explicit constraints $F$ may be extended to a mapping from $\R^n \times \R^m$ to $\R^n \times \R^m$
with $(x,z) \mapsto F(x,z)$ and $z \in \R^m$ acting as a multiplier for the constraints.
Assuming that $F(x,z)=0$ describes a KKT-system (typically implying that $F$ is not differentiable), we identify for a given $\mu > 0$ the equation $F(x,z;\mu)=0$ with the
\textit{perturbed} KKT-system (preserving differentiability) as defined in \eqref{eq:pd_system}.
For the model problem \eqref{eq_abstract_discr} this corresponds to $F_{box}(x,z^a,z^b;\mu)=0$.

Consequently, we may define the (primal-dual) barrier homotopy map as
\begin{align}
  H(x,z^a,z^b,\mu,t) &:= F_{box}(x,z^a,z^b;\mu)- (1-t)F_{box}(x_0,z_0^a,z_0^b;\mu_0) \label{def:barrier_hmtpy} \\
  &=
  \begin{bmatrix}
      \nabla f(x) - z^a + z^b - (1-t) (\nabla f(x_0) - z_0^a + z_0^b) \\
       z^a \odot c^a(x) - \mu e_n - (1-t) (z_0^a \odot c^a(x_0) - \mu_0 e_n) \\
      z^b \odot c^b(x) - \mu e_n - (1-t) (z_0^b \odot c^b(x_0) - \mu_0 e_n)
  \end{bmatrix} \nonumber\\
  &=
  \begin{bmatrix}
      \nabla f(x) - z^a + z^b - (1-t) (\nabla f(x_0) - z_0^a + z_0^b) \\
       z^a \odot c^a(x) - \mu e_n \\
      z^b \odot c^b(x) - \mu e_n
  \end{bmatrix}
  \nonumber
\end{align}
where $z_0^a \odot c^a(x_0) - \mu_0 e_n = z_0^b \odot c^b(x_0) - \mu_0 e_n = 0$ holds by construction of $z_0^a$ and $z_0^b$.
The Jacobian of $H$ with respect to the homotopy variable $(x,z^a,z^b)$ is further given by
\begin{align}
   JH(x,z^a,z^b,\mu,t)=
  \begin{bmatrix}
    \nabla^2 f(x) & -I_n &I_n \\
    Z^a & C^a(x) & 0 \\
    -Z^b & 0 & C^b(x)
  \end{bmatrix}
\end{align}
with the diagonal matrices $Z^a,Z^b$ and $C^a(x),C^b(x)$ defined as in \eqref{eq:pd_newton}. The Newton system thus coincides with the one in the primal-dual barrier approach \eqref{eqn:pd_newton}
with the only difference of the perturbation of the right-hand side controlled by $t$ occurring in the homotopy approach exclusively.

\begin{remark}
  \label{rem:barrier_rule}
  The update rule for the barrier parameter $\mu$ is designed to depend on the homotopy parameter $t$, i.e., $\mu = \mu(t)$.
  Correspondingly, it is required that $\mu(0)=\mu_0$ and $\mu(1)=\mu_\infty$ with $\mu_0$ the initial and $\mu_\infty$ the target barrier prameters, e.g.,
  $\mu(t):= t \mu_\infty + (1-t) \mu_0$.
  The homotopy map in \eqref{def:barrier_hmtpy} can thus be expressed as $H(x,z^a,z^b,t):= H(x,z^a,z^b,\mu(t),t)$ eliminating the explicit dependence on $\mu$.
  Therefore by providing the mapping $t \mapsto \mu(t)$ and replacing $x$ by $(x,z^a,z^b)$ \Cref{alg:pc_hmtpy} can be employed to compute a
  stationary point of \eqref{eq_abstract_discr}.
\end{remark}

\paragraph{Extension to PDE constraints}

Neglecting the box constraints on the design variable the Lagrangian $\Lc:Q \times V \times V \to \R$ of the PDE-constrained minimization problem \eqref{pr:red_abstract} is given by
$\Lc(\rho,u,p):= J(\rho,u) + \innerpr{e(\rho,u)}{p}_{V^* \times V}$.
With $\uprho, \mathrm{u}$ and $\mathrm{p}$ denoting the respective finite element coefficient vectors of $\rho,u$ and $p$ the stationarity condition of
the discretized Lagrangian $\mathrm{L}:\R^n \times \R^l \times \R^l \to \R$ is then given by
\begin{align}
   \gradL:=
  \begin{bmatrix}
    \dLn{\uprho} \\
    \dLn{u} \\
    \dLn{p}
  \end{bmatrix}
  =0
\end{align}
where $\dLn{p}=0$ and $\dLn{u}=0$ represent the discretized state and adjoint equation.

We are now able to implicitly enforce the box constraints on $\uprho$ by embedding $F(\uprho,\mathrm{u},\mathrm{p}):=\gradL$ into the primal-dual
barrier framework, i.e.,
\begin{align}
  F_{box}(\uprho,\mathrm{u},\mathrm{p},z^a,z^b;\mu) :=
  \begin{bmatrix}
     \dLn{\uprho} - z^a + z^b\\
     \dLn{u} \\
     \dLn{p} \\
     z^a \odot c^a(\uprho) - \mu e_n \\
    z^b \odot c^b(\uprho) - \mu e_n
  \end{bmatrix}.
\end{align}

The corresponding global homotopy is further defined as
\begin{align}
  H(\uprho,\mathrm{u},\mathrm{p},z^a,z^b,t) &= H(\uprho,\mathrm{u},\mathrm{p},z^a,z^b,\mu(t),t)\\
   :&= F_{box}(\uprho,\mathrm{u},\mathrm{p},z^a,z^b;\mu(t)) - (1-t) F_{box}(\uprho_0,\mathrm{u}_0,\mathrm{p}_0,z_0^a,z_0^b;\mu_0) \nonumber\\
   &= \begin{bmatrix}
    \dLn{\uprho} - z^a + z^b - (1-t)(\dLninit{\uprho} - z_0^a + z_0^b)\\
    \dLn{u} - (1-t)(\dLninit{u})\\
    \dLn{p} - (1-t)(\dLninit{p})\\
     z^a \odot c^a(\uprho) - \mu(t) e_n - (1-t)(z_0^a \odot c^a(\uprho_0) - \mu_0 e_n)\\
      z^b \odot c^b(\uprho) - \mu(t) e_n - (1-t)( z_0^b \odot c^b(\uprho_0) - \mu_0 e_n)
   \end{bmatrix} \nonumber\\
   &=
   \begin{bmatrix}
    \dLn{\uprho} - z^a + z^b - (1-t)(\dLninit{\uprho} - z_0^a + z_0^b)\\
    \dLn{u} \\
    \dLn{p} \\
    z^a \odot c^a(\uprho) - \mu(t) e_n \\
    z^b \odot c^b(\uprho) - \mu(t) e_n
   \end{bmatrix} \nonumber
\end{align}
where we additionally used $\dLninit{p}=0$ and $\dLninit{u}=0$, i.e., the initializations of the state and adjoint variable satisfy the respective equations.

The Jacobian of the homotopy map with respect to $(\uprho,\mathrm{u},\mathrm{p},z^a,z^b)$ required for the Newton system can then be stated as
\begin{align}
   JH(\uprho,\mathrm{u},\mathrm{p},z^a,z^b,t) =
   \begin{bmatrix}
    \ddLn{\uprho}{\uprho}  & \ddLn{u}{\uprho}  & \ddLn{p}{\uprho} & -I_n & I_n \\
    \ddLn{\uprho}{u}       & \ddLn{u}{u} & \ddLn{p}{u} & 0 & 0 \\
    \ddLn{\uprho}{p}      & \ddLn{u}{p}  & 0 & 0 & 0  \\
    Z^a & 0 & 0 & C^a(\uprho) & 0 \\
    -Z^b & 0 & 0 & 0& C^b(\uprho)
  \end{bmatrix}
   \label{eq:newton_barrier_hmtpy_pde}
\end{align}
where we used that $\ddLn{p}{p}=0$.

Consequently, with the initial vectors $\mathrm{u}_0$ and $\mathrm{p}_0$ computed based on $\uprho_0$, with the initialization of the dual variables $z^a$ and $z^b$ as
described in \Cref{sec:barrier} and with an update rule on the barrier parameter, an extension of \Cref{alg:pc_hmtpy} can be employed to compute a stationary point of the
original target problem \eqref{pr:red_abstract}. This procedure is summarized in \Cref{alg:barrier_hmtpy_pde}.

\begin{figure}
  \centering
  \begin{minipage}{0.75\linewidth}
    \begin{algorithm}[H]
      \caption{Barrier homotopy method for a PDE constrained optimization problem}\label{alg:barrier_hmtpy_pde}
      \begin{algorithmic}
        \Require $\uprho_0 \in (0,1)^n$, update rules $\theta(t)$, $\mu(t)$, $\mu_0 \gg 0$, $\mu_\infty \approx 0$
        \Ensure $\begin{pmatrix}
          \uprho, \mathrm{u}, \mathrm{p}, z^a, z^b
        \end{pmatrix}$ approximately fulfilling the KKT-system of the discretization of \eqref{pr:red_abstract}
        \State $\uprho \gets \uprho_0$, $t \gets 0$, $\mu \gets \mu_0 =: \mu(0)$
      \State $\mathrm{u}_0 \gets \text{solution of }  \partial_{\mathrm{p}} \mathrm{L}(\uprho_0,\mathrm{u},0)=0$, $\mathrm{p}_0 \gets \text{solution of }  \partial_{\mathrm{p}} \mathrm{L}(\uprho_0,\mathrm{u}_0,\mathrm{p})=0$
      \State $z^a \gets \mu_0 \oslash c^a(\uprho_0)$, $z^b \gets \mu_0 \oslash c^b(\uprho_0)$
      \While{$t < 1$}
      \State $t \gets \theta(t)$, $\mu \gets \mu(t)$
      \While{$\Norm{H(\uprho,\mathrm{u},\mathrm{p},z^a,z^b,t)}>tol$}
      \State Solve $JH(\uprho,\mathrm{u},\mathrm{p},z^a,z^b,t)
      \begin{bmatrix}
        \Delta \uprho \\
        \Delta \mathrm{u} \\
        \Delta \mathrm{p} \\
        \Delta z^a  \\
        \Delta z^b
      \end{bmatrix} = -H(\uprho,\mathrm{u},\mathrm{p},z^a,z^b,t)$
      \State $\begin{bmatrix}
        \uprho \\
        \mathrm{u} \\
        \mathrm{p} \\
        z^a  \\
        z^b
      \end{bmatrix} \gets
      \begin{bmatrix}
        \uprho \\
        \mathrm{u} \\
        \mathrm{p} \\
        z^a  \\
        z^b
      \end{bmatrix} + \begin{bmatrix}
        \Delta \uprho \\
        \Delta \mathrm{u} \\
        \Delta \mathrm{p} \\
        \Delta z^a  \\
        \Delta z^b
      \end{bmatrix}$
      \EndWhile
      \EndWhile
    \end{algorithmic}
  \end{algorithm}
  \end{minipage}
\end{figure}

\section{Numerical Results}{\label{sec:results}}

For the subsequent numerical tests we intend to compute a stationary point of the compliance problem from \Cref{sec:model_prb} as defined in \eqref{pr:red_abstract} by applying
the barrier homotopy approach specified in \Cref{alg:barrier_hmtpy_pde} to its discretized formulation. Numerical experiments have shown that in this setting the overall number of required homotopy steps does not decrease significantly when emplyoing a first order predictor (which entails additional computational effort for solving the equation \eqref{eq:h_predictor}) instead of a zero order predictor (see \Cref{rem:h_predictor}). We therefore restrict ourselves to the utilization of the latter subsequently.



\begin{wrapfigure}{r}{6.5cm}
					\begin{tikzpicture}[scale=3.2]
						\draw[-] (0,0) -- (2,0);
						\draw[-] (0,0) -- (0,1);
						\draw[-] (2,0) -- (2,1);
						\draw[-] (0,1) -- (2,1);

						\draw[-] (0.2,-0.05) -- (0.2,0.05);

						\draw[-] (0.00,-0.03) -- (0.02,0.03);
						\draw[-] (0.02,-0.03) -- (0.04,0.03);
						\draw[-] (0.04,-0.03) -- (0.06,0.03);
						\draw[-] (0.06,-0.03) -- (0.08,0.03);
						\draw[-] (0.08,-0.03) -- (0.10,0.03);
						\draw[-] (0.10,-0.03) -- (0.12,0.03);
						\draw[-] (0.12,-0.03) -- (0.14,0.03);
						\draw[-] (0.14,-0.03) -- (0.16,0.03);
						\draw[-] (0.16,-0.03) -- (0.18,0.03);
						\draw[-] (0.18,-0.03) -- (0.20,0.03);

						\draw[-] (0.9,-0.05) -- (0.9,0.05);

						\draw[-] (1.1,-0.05) -- (1.1,0.05);

						\draw[-] (1.8 + 0.00,-0.03) -- (1.8 + 0.02,0.03);
						\draw[-] (1.8 + 0.02,-0.03) -- (1.8 + 0.04,0.03);
						\draw[-] (1.8 + 0.04,-0.03) -- (1.8 + 0.06,0.03);
						\draw[-] (1.8 + 0.06,-0.03) -- (1.8 + 0.08,0.03);
						\draw[-] (1.8 + 0.08,-0.03) -- (1.8 + 0.10,0.03);
						\draw[-] (1.8 + 0.10,-0.03) -- (1.8 + 0.12,0.03);
						\draw[-] (1.8 + 0.12,-0.03) -- (1.8 + 0.14,0.03);
						\draw[-] (1.8 + 0.14,-0.03) -- (1.8 + 0.16,0.03);
						\draw[-] (1.8 + 0.16,-0.03) -- (1.8 + 0.18,0.03);
						\draw[-] (1.8 + 0.18,-0.03) -- (1.8 + 0.20,0.03);

						\draw[-] (1.8,-0.05) -- (1.8,0.05);
						\draw (1.0,0.05) node[above,scale=0.55] {$\Gamma_{N,g_N}$};
						\draw (0.1,0.05) node[above,scale=0.55] {$\Gamma_{D,0}$};
						\draw (1.9,0.05) node[above,scale=0.55] {$\Gamma_{D,0}$};

						\draw[thick,->] (1.1,-0.08) -- (1.1,-0.3) node[right,pos=0.5,scale=0.55] {$g_N$};
            \draw[thick,->] (1.0,-0.08) -- (1.0,-0.3);
            \draw[thick,->] (0.9,-0.08) -- (0.9,-0.3);


					\end{tikzpicture}
				\caption{The design domain $\Omega=[0,2.4] \times [0,0.8]$}
        \label{fig:omega}
\end{wrapfigure}
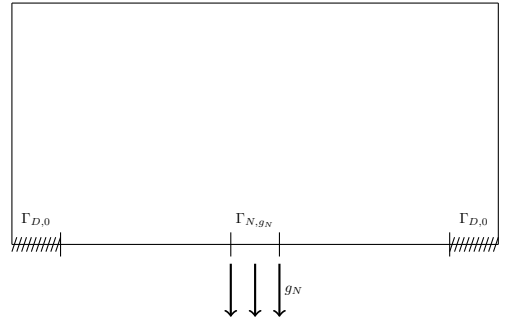

The boundary of the rectangular computational domain $\Omega=[0,2.4] \times [0,0.8]$ (see \Cref{fig:omega}) consists of the homogeneous Dirichlet boundary $\Gamma_{D,0}=[0,0.12] \times \{0\} \cup [2.28,2.4] \times \{0\}$,
the nonhomogeneous Neumann boundary $\Gamma_{N,g_N}= [1.08,1.32] \times \{0\}$ where the traction force ${g_N = (0,-1)^\top}$ is applied and the homogeneous Neumann boundary on
the remaining parts $\partial \Omega \setminus (\Gamma_{D,0} \cup \Gamma_{N,g_N})$.



For the discretization of $\Omega$ a mesh with $11100$ triangular elements and $5711$ vertices is used.
The function space $Q=H^1(\Omega)$ for the density variable and the vector valued space $V=(H_{0,\Gamma_D}^1(\Omega))^2$
for the state and adjoint variable are both approximated by piecewise linear and globally continuous functions. The Lamé parameters of the first material
and the ersatz material (imitating air) are defined as $\lambda^L_1:=\num{0.750}$, $\mu^L_1:=\num{0.375}$ and $\lambda^L_0:=\num{7.498e-5}$, $\mu^L_0:=\num{3.750e-5}$,
respectively.

The density variable representing the initial (nonbinary) design is set to the constant $\uprho_0 \equiv 0.5$. To enforce a non-full design
the weight of the volume control term is set to $\gamma=9.75$. The Ginzburg-Landau parameter (which is proportional to the thickness of the interface between the two materials)
is defined as $\varepsilon=\num{0.0075}$ and the weight of the corresponding energy term is set to $\beta=0.5$. For the barrier parameter $\mu$
an interpolation linear in $t$ between the initial value $\mu_0=\num{50}$ and the target value $\mu_\infty=\num{0.001}$ is chosen. The update of the
homotopy parameter obeys a naive stepsize rule with the initial step $\Delta t_{init}=0.25$. In case of convergence of Newton's method after the
update on $t$, the stepsize is increased by $50\%$ (restricted by the upper bound $\Delta t_{max}=0.25$); if Newton's method diverges, the stepsize is halvened.

\begin{figure}

\newcommand{\myscale}{0.16}

\begin{figure}[H]
\begin{minipage}{0.33\textwidth}
  \begin{figure}[H]
	\centering
  \captionsetup{justification=centering,font=footnotesize}
	\includegraphics[scale=\myscale, trim={0.5cm 8.3cm 0.5cm 8.3cm},clip]{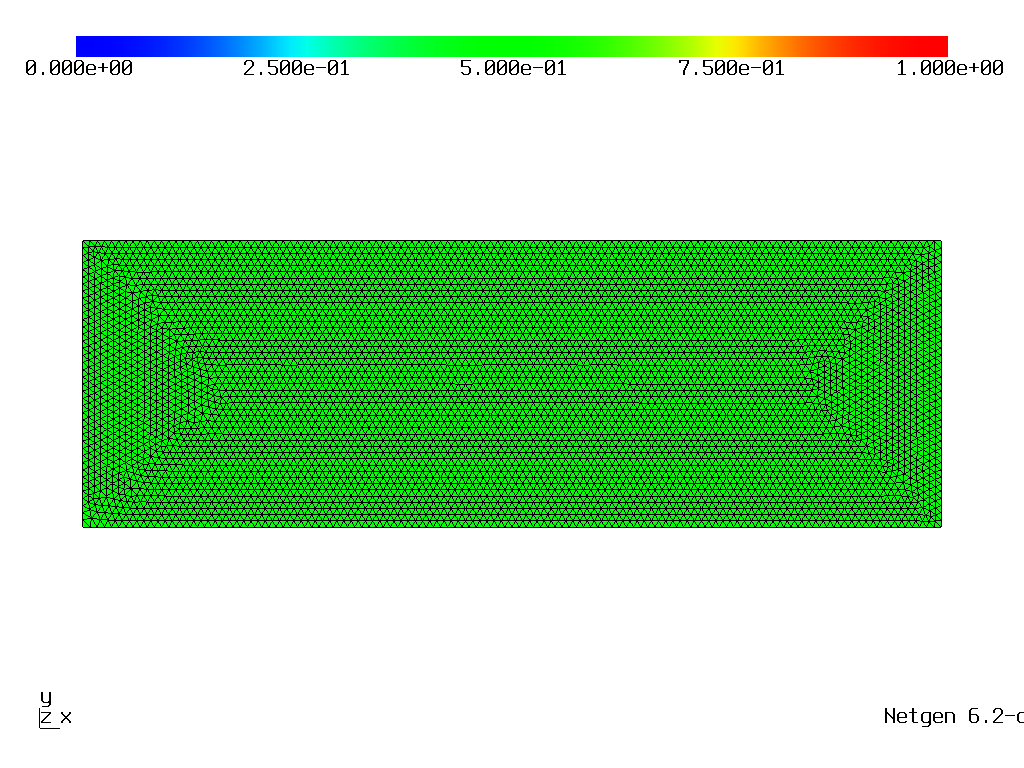}
  \caption*{$t=0.000000$}
	\end{figure}
\end{minipage}
\hfil
\begin{minipage}{0.33\textwidth}
  \begin{figure}[H]
	\centering
  \captionsetup{justification=centering,font=footnotesize}
	\includegraphics[scale=\myscale, trim={0.5cm 8.3cm 0.5cm 8.3cm},clip]{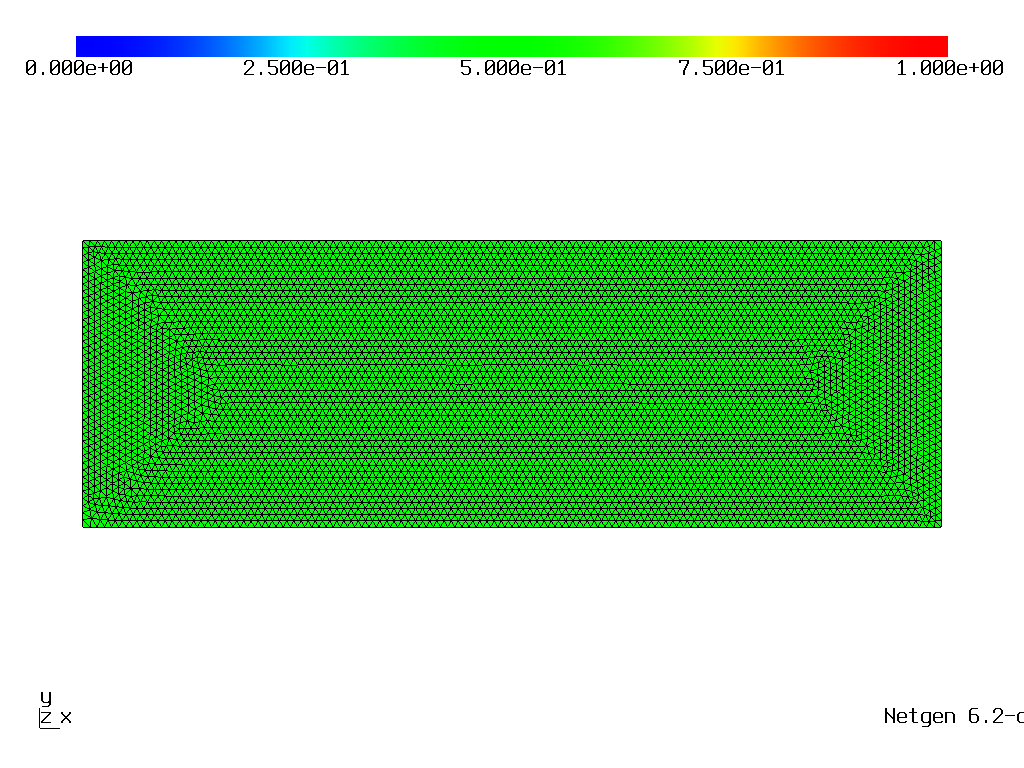}
  \caption*{$t=0.500000$}
	\end{figure}
\end{minipage}
\hfil
\begin{minipage}{0.33\textwidth}
  \begin{figure}[H]
	\centering
  \captionsetup{justification=centering,font=footnotesize}
	\includegraphics[scale=\myscale, trim={0.5cm 8.3cm 0.5cm 8.3cm},clip]{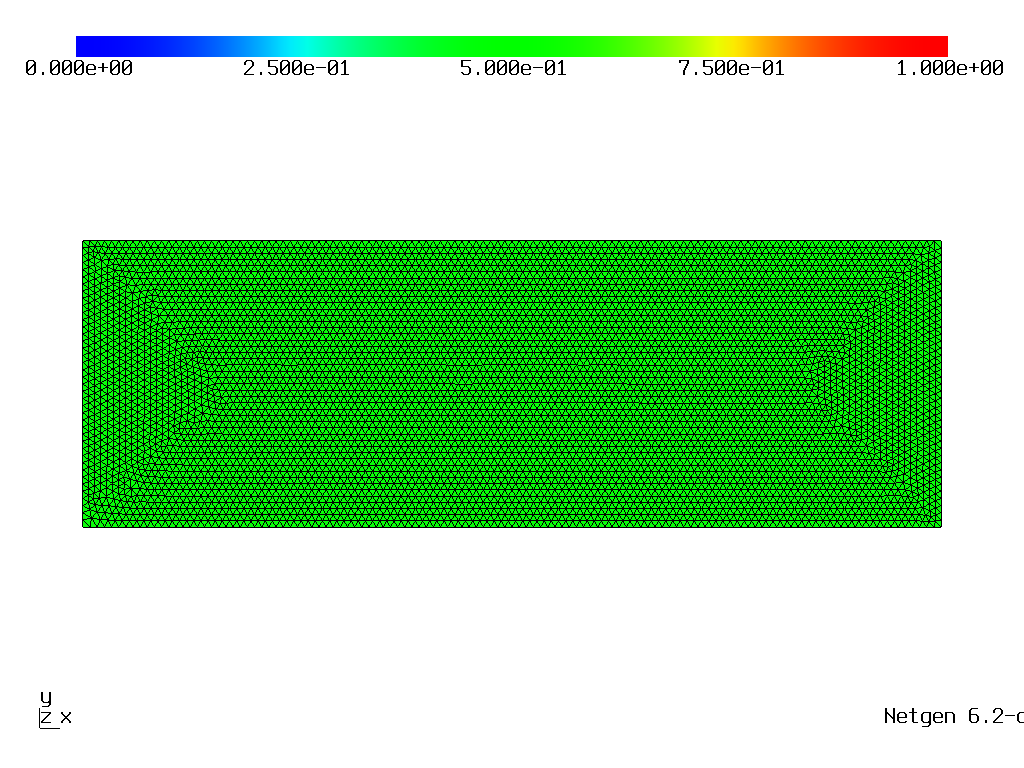}
  \caption*{$t=0.937500$}
	\end{figure}
\end{minipage}
\end{figure}

\vspace{-0.5cm}

\begin{figure}[H]
\begin{minipage}{0.33\textwidth}
  \begin{figure}[H]
	\centering
  \captionsetup{justification=centering,font=footnotesize}
	\includegraphics[scale=\myscale, trim={0.5cm 8.3cm 0.5cm 8.3cm},clip]{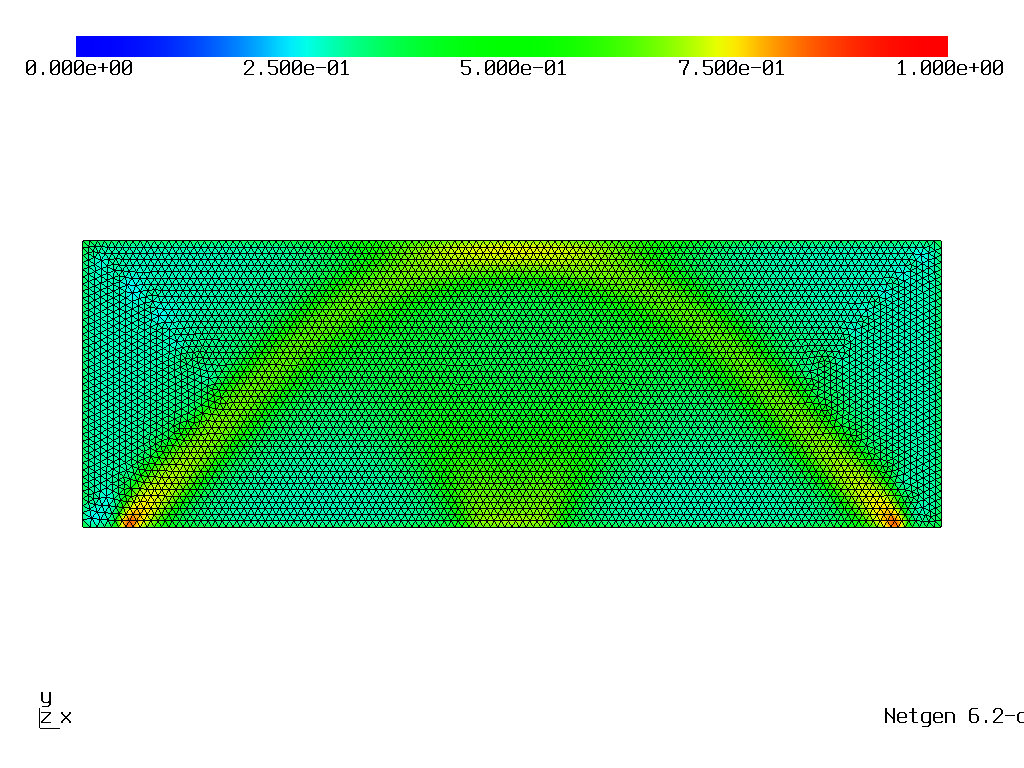}
  \caption*{$t=0.999931$}
	\end{figure}
\end{minipage}
\hfil
\begin{minipage}{0.33\textwidth}
  \begin{figure}[H]
	\centering
  \captionsetup{justification=centering,font=footnotesize}
	\includegraphics[scale=\myscale, trim={0.5cm 8.3cm 0.5cm 8.3cm},clip]{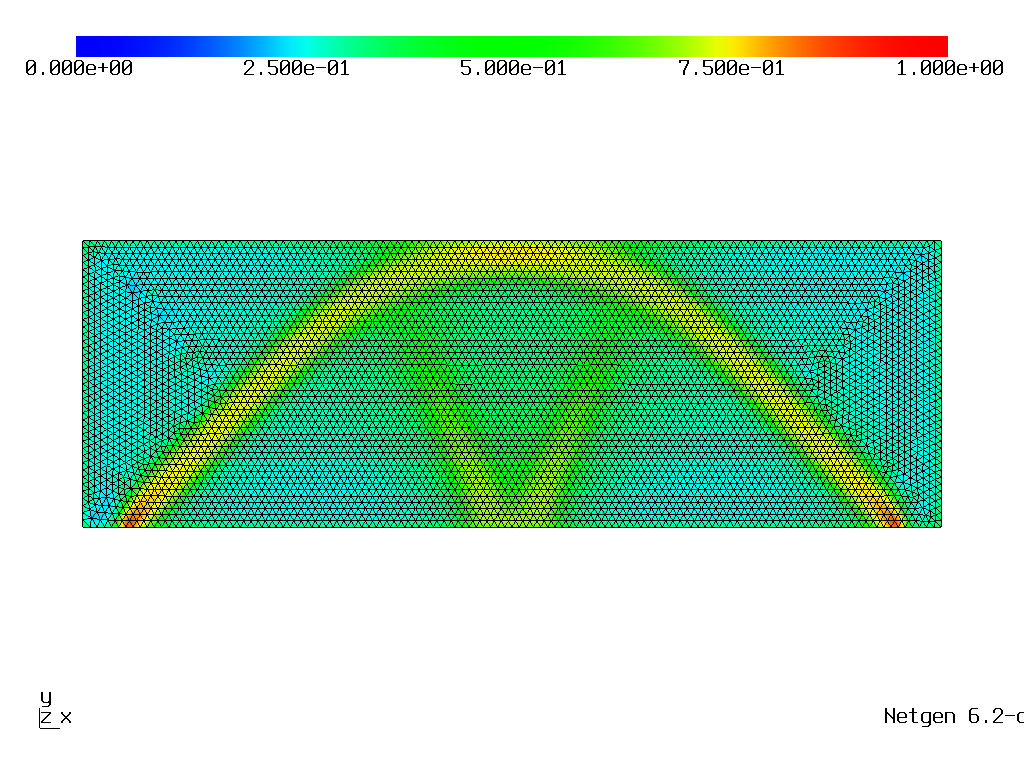}
  \caption*{$t=0.999946$}
	\end{figure}
\end{minipage}
\hfil
\begin{minipage}{0.33\textwidth}
  \begin{figure}[H]
	\centering
  \captionsetup{justification=centering,font=footnotesize}
	\includegraphics[scale=\myscale, trim={0.5cm 8.3cm 0.5cm 8.3cm},clip]{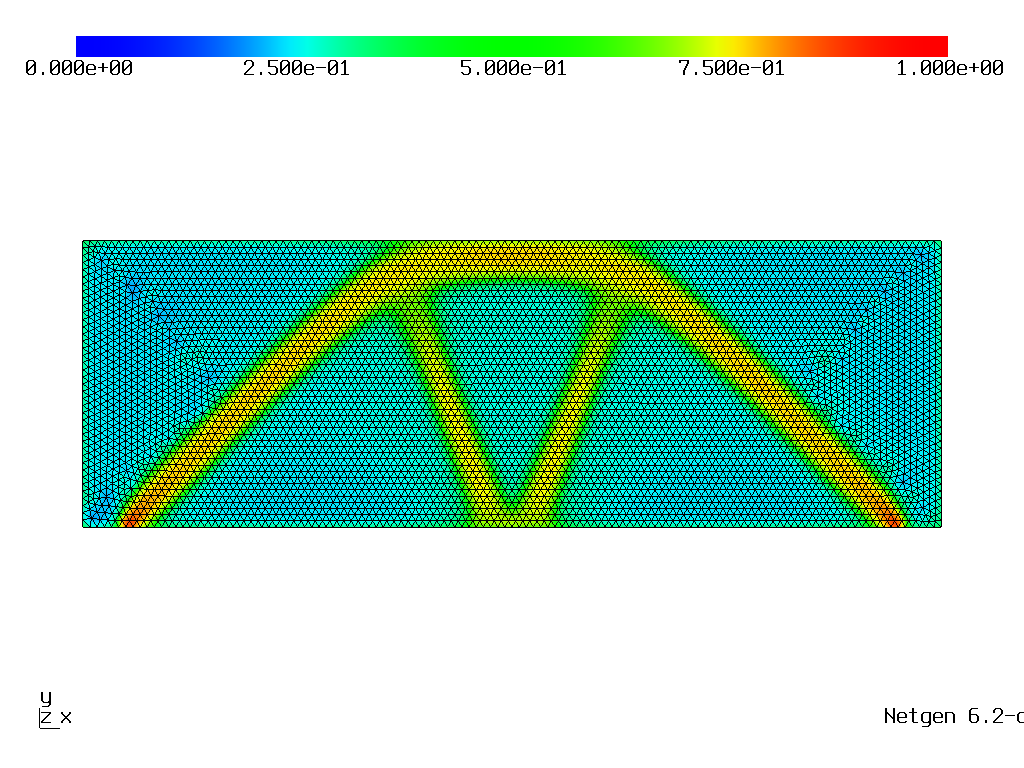}
  \caption*{$t=0.999956$}
	\end{figure}
\end{minipage}
\end{figure}

\vspace{-0.5cm}

\begin{figure}[H]
  \captionsetup{justification=centering}
\begin{minipage}{0.33\textwidth}
  \begin{figure}[H]
	\centering
  \captionsetup{justification=centering,font=footnotesize}
	\includegraphics[scale=\myscale, trim={0.5cm 8.3cm 0.5cm 8.3cm},clip]{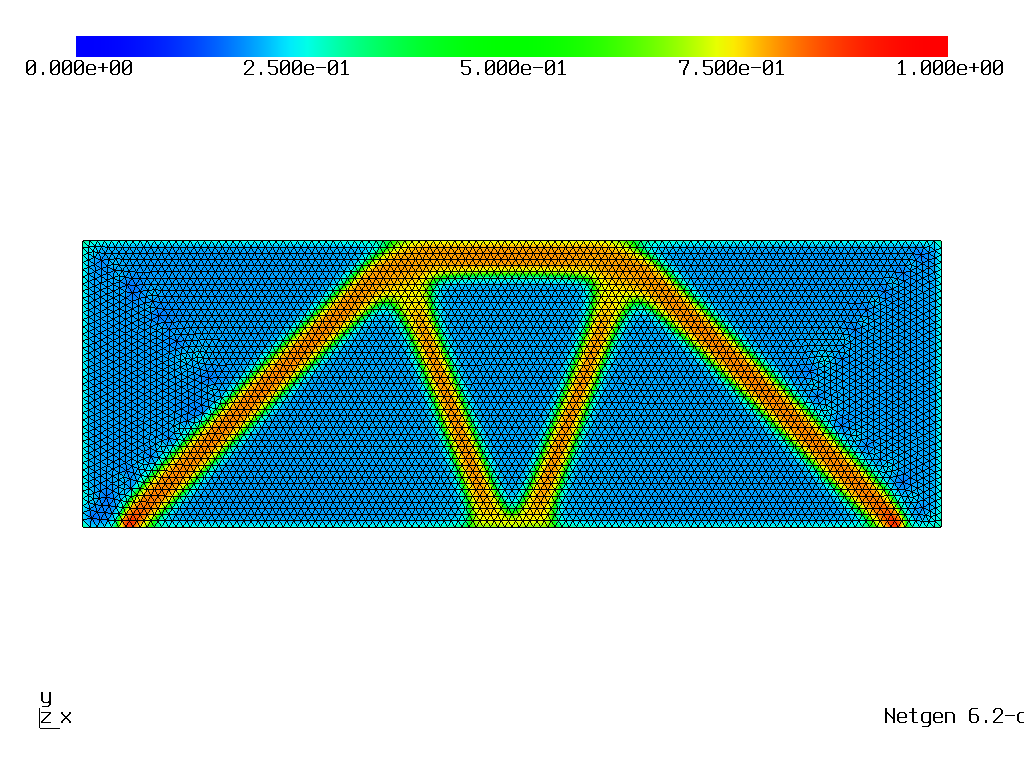}
  \caption*{$t=0.999974$}
	\end{figure}
\end{minipage}
\hfil
\begin{minipage}{0.33\textwidth}
  \begin{figure}[H]
	\centering
	\includegraphics[scale=\myscale, trim={0.5cm 8.3cm 0.5cm 8.3cm},clip]{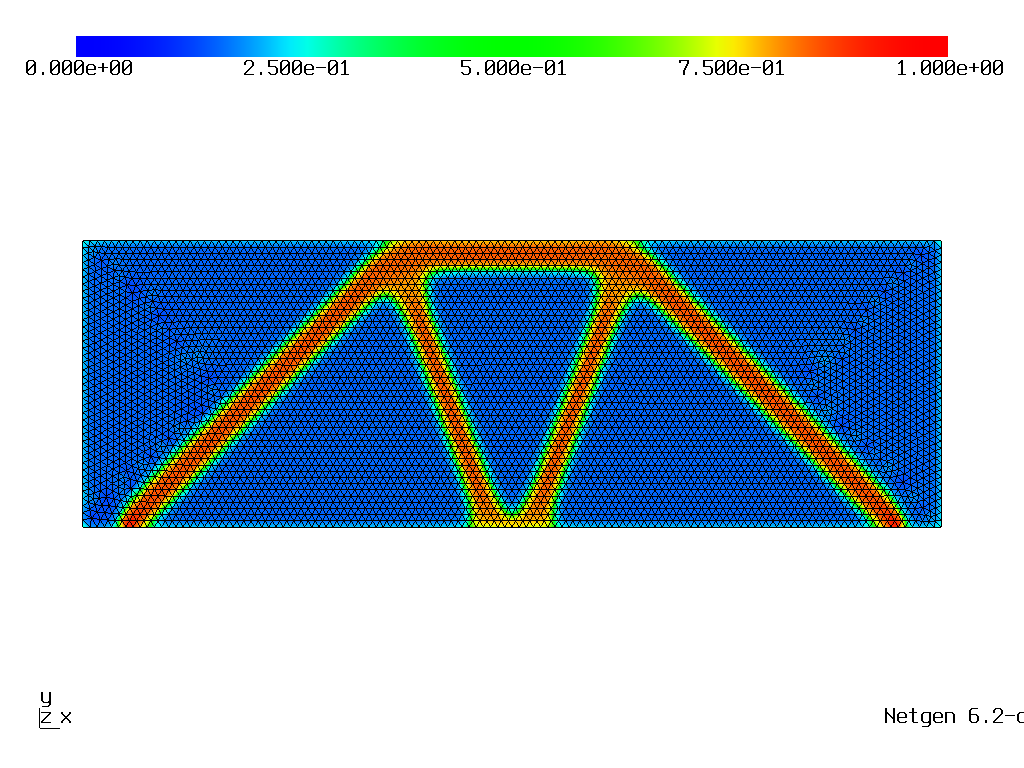}
  \captionsetup{justification=centering,font=footnotesize}
  \caption*{$t=0.999988$}
	\end{figure}
\end{minipage}
\hfil
\begin{minipage}{0.33\textwidth}
  \begin{figure}[H]
	\centering
	\includegraphics[scale=\myscale, trim={0.5cm 8.3cm 0.5cm 8.3cm},clip]{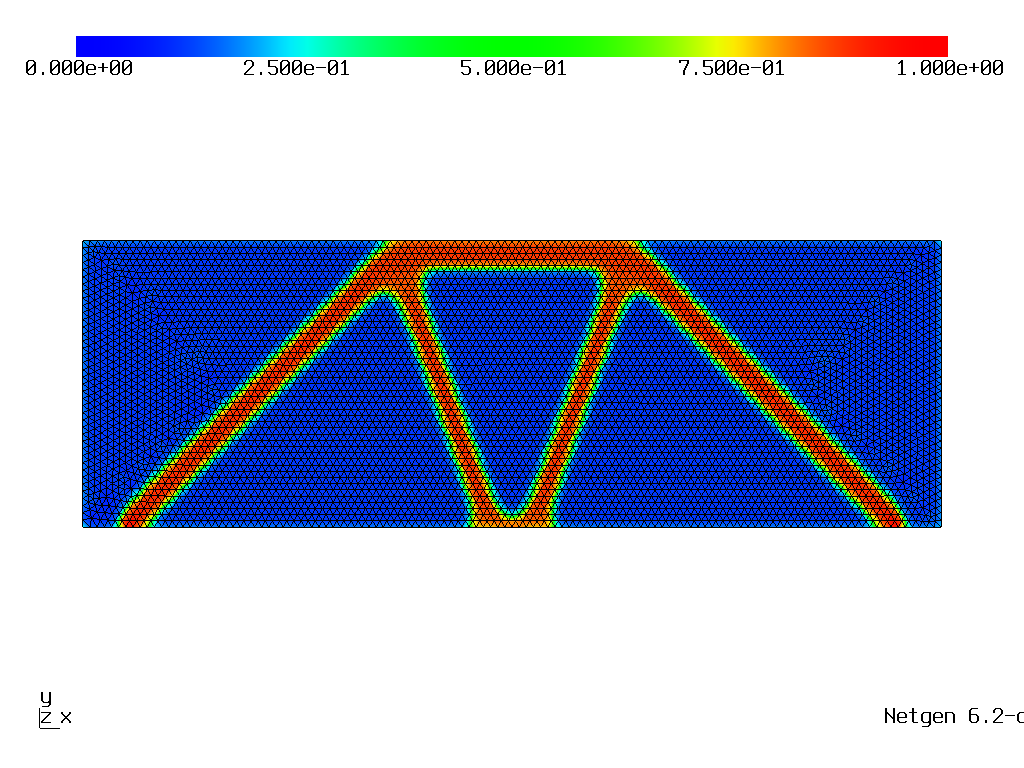}
  \captionsetup{justification=centering,font=footnotesize}
  \caption*{$t=1.000000$}
	\end{figure}
\end{minipage}
\end{figure}

\caption{Evolution of the design starting out from an intermediate material ($\uprho_0 \equiv 0.5$) in green and approximating a binary design with the first material in red and the ersatz material in blue as $t$ approaches the target parameter $1$}
\label{fig:evol_design}
\end{figure}

The method defined as above required $26$ (successful) iterations, i.e., steps in $t$, to compute a stationary point.
Owing to the fact that by the naive update rule the stepsize was sometimes set too large, in total $47$ iterations were conducted (see \Cref{fig:t_mu_history}).
The evolution of the associated designs (represented by the values of the density variable) depending on $t$ is depicted in \Cref{fig:evol_design} where
red corresponds to the first material, blue to the ersatz material and the colours in between to an intermediate material.
It is apparent in \Cref{fig:evol_design} and \Cref{fig:t_mu_history} that the small stepsize necessary to cause Newton's method to converge as $t$ approaches $1$ correlates
with a significant change in the design near $t=1$. For smaller values of $t$ this pattern cannot be observed.

\begin{figure}
  \begin{minipage}{0.48\textwidth}
  \begin{figure}[H]
  \centering
  \captionsetup{justification=centering}
  \begin{tikzpicture}
  \begin{axis}[xlabel={iteration index},
               ylabel={$t$},
               xtick={0,10,20,30,40},
               minor xtick = {1,2,3,4,5,6,7,8,9,11,12,13,14,15,16,17,18,19,21,22,23,24,25,26,27,28,29,31,32,33,34,35,36,37,38,39,41,42,43,44,45,46,47},
               minor ytick={0.05,0.1,0.15,0.25,0.3,0.35,0.45,0.5,0.55,0.65,0.7,0.75,0.85,0.9,0.95}]
    \addplot [black,-, mark=*, mark size=1pt] table [x=it, y=t, col sep=comma] {param_history.csv};
  \end{axis}
\end{tikzpicture}
\end{figure}
\end{minipage}
\hfil
\begin{minipage}{0.48\textwidth}
  \begin{figure}[H]
  \centering
  \captionsetup{justification=centering}
  \begin{tikzpicture}
  \begin{axis}[xlabel={iteration index},
               ylabel={$\mu$},
               xtick={0,10,20,30,40},
               minor xtick = {1,2,3,4,5,6,7,8,9,11,12,13,14,15,16,17,18,19,21,22,23,24,25,26,27,28,29,31,32,33,34,35,36,37,38,39,41,42,43,44,45,46,47}]
    \addplot [black,-, mark=*, mark size=1pt] table [x=it, y=mu, col sep=comma] {param_history.csv};
  \end{axis}
\end{tikzpicture}
\end{figure}
\end{minipage}
\captionsetup{justification=centering}
\caption{Evolution of the homotopy parameter (left) and the barrier parameter (right)}
\label{fig:t_mu_history}
\end{figure}

\section{Conclusion and Outlook}

We developed a (primal-dual) barrier homotopy approach and employed it to compute a stationary point of a nonconvex design optimization problem with a PDE constraint acting on the state variable and box constraints restricting the admissible values of the density variable. This continuation approach involves a twofold perturbation of the target problem. The first one is related to an auxiliary problem which is induced by the choice of an initial guess and the
second one corresponds to a barrier introduced to prevent the design variable from escaping the feasible region described by the box constraints.
In the course of the method the weight of the barrier and the contribution of the auxiliary problem is simultaneously decreased (controlled by the homotopy parameter) until the target problem is reached.

Numerical results for a compliance minimization problem validate
the utilization of both the homotopy and the barrier approach for its globalization property (in terms of the initial guess) and its implicit handling of the inequality constraints, respectively.
That is, starting with an unbiased nonbinary design (and therefore not requiring any a priori knowledge about a locally optimal design) leads to the detection of a feasible stationary solution
whose topology is not already inherent in the initial guess.

It is worth noting that for the choice of the target barrier parameter $\mu_\infty$ the balance between being sufficiently small to approximate the target problem well and
being large enough to facilitate the computation of solutions of the subproblems as $\mu$ approaches $\mu_\infty$ is important for the setup of the algorithm.
This is observed in numerical tests where decreasing $\mu_\infty$ below a certain threshold does have the effect of allowing the values of the density variable closer to $0$ or $1$
but not inducing topological changes anymore and amplifying the ill-conditioning of the target problem instead.

In future work, we plan to employ the presented approach for efficiently exploring the design space of multi-objective topology optimization problems. Combining the tracing of a Pareto curve \cite{schmidt2008pareto,MartinSchuetze2018} with deflation \cite{PapadopoulosFarrellSurowiec2021} may yield better locally optimal Pareto fronts compared to more conventional first order approaches.

\section*{Acknowledgment}
The work of the authors is partially supported by the joint DFG/FWF Collaborative Research Centre CREATOR (DFG: Project-ID 492661287/TRR 361; FWF: 10.55776/F90) at TU Darmstadt, TU Graz, JKU Linz and RICAM Linz. P.G. is partially supported by the State of Upper Austria.

\bibliography{arxiv.bib}

@Book{BendsoeSigmund2003,
 title = {Topology Optimization: Theory, Methods and Applications},
 publisher = {Springer, Berlin},
 year = {2003},
 author = {M. P. Bends\o{}e and O. Sigmund},
}

@article{allaire2020survey,
  title={Shape and topology optimization},
  author={Allaire, G. and Dapogny, C. and Jouve, F.},
  journal={in Geometric partial differential equations, part II, A. Bonito and R. Nochetto eds., Handbook of Numerical Analysis, vol. 22},
  pages = {1--132},
  year={2021}
 }

@article{allaire2004structural,
  title={Structural optimization using sensitivity analysis and a level-set method},
  author={Allaire, G. and Jouve, F. and Toader, A.-M.},
  journal={Journal of computational physics},
  volume={194},
  number={1},
  pages={363--393},
  year={2004},
  publisher={Elsevier}
}

@article{Cherrire2022,
  title = {Multi-material topology optimization using Wachspress interpolations for designing a 3-phase electrical machine stator},
  volume = {65},
  ISSN = {1615-1488},
  url = {http://dx.doi.org/10.1007/s00158-022-03460-1},
  DOI = {10.1007/s00158-022-03460-1},
  number = {12},
  journal = {Structural and Multidisciplinary Optimization},
  publisher = {Springer Science and Business Media LLC},
  author = {Cherrière,  T. and Laurent,  L. and Hlioui,  S. and Louf,  F. and Duysinx,  P. and Geuzaine,  C. and Ben Ahmed,  H. and Gabsi,  M. and Fernández,  E.},
  year = {2022},
  month = nov
}

@misc{GanglKrennDeGersem2025,
      title={Multi-material topology optimization of electric machines under maximum temperature and stress constraints},
      author={Gangl, P. and Krenn, N. and De Gersem, H.},
      year={2025},
      eprint={2504.12426},
      archivePrefix={arXiv},
      primaryClass={math.OC},
      url={https://arxiv.org/abs/2504.12426},
}

@article{FepponAllaireDapognyJolivet2020,
title = {Topology optimization of thermal fluid–structure systems using body-fitted meshes and parallel computing},
journal = {Journal of Computational Physics},
volume = {417},
pages = {109574},
year = {2020},
issn = {0021-9991},
doi = {https://doi.org/10.1016/j.jcp.2020.109574},
url = {https://www.sciencedirect.com/science/article/pii/S002199912030348X},
author = {F. Feppon and G. Allaire and C. Dapogny and P. Jolivet}
}

@Article{BendsoeSigmund1999,
	author="Bends{\o}e, M. P.
	and Sigmund, O.",
	title="Material interpolation schemes in topology optimization",
	journal="Archive of Applied Mechanics",
	year="1999",
	volume="69",
	number="9",
	pages="635--654",
	issn="1432-0681",
	doi="10.1007/s004190050248",
	url="http://dx.doi.org/10.1007/s004190050248"
}

@article{PapadopoulosFarrellSurowiec2021,
author = {Papadopoulos, I. P. A. and Farrell, P. E. and Surowiec, T. M.},
title = {Computing Multiple Solutions of Topology Optimization Problems},
journal = {SIAM Journal on Scientific Computing},
volume = {43},
number = {3},
pages = {A1555-A1582},
year = {2021},
doi = {10.1137/20M1326209},
URL = { https://doi.org/10.1137/20M1326209},
eprint = { https://doi.org/10.1137/20M1326209}
}

@article {Svanberg1987,
author = {Svanberg, K.},
title = {The method of moving asymptotes -- a new method for structural optimization},
journal = {International Journal for Numerical Methods in Engineering},
volume = {24},
number = {2},
publisher = {John Wiley & Sons, Ltd},
issn = {1097-0207},
url = {http://dx.doi.org/10.1002/nme.1620240207},
doi = {10.1002/nme.1620240207},
pages = {359--373},
year = {1987},
}

@article{AdamSurowiec2018,
  title = {A {PDE}-constrained optimization approach for topology optimization of strained photonic devices},
  volume = {19},
  ISSN = {1573-2924},
  url = {http://dx.doi.org/10.1007/s11081-018-9394-5},
  DOI = {10.1007/s11081-018-9394-5},
  number = {3},
  journal = {Optimization and Engineering},
  publisher = {Springer Science and Business Media LLC},
  author = {Adam,  L. and Hinterm\"{u}ller,  M. and Surowiec,  T. M.},
  year = {2018},
  month = jul,
  pages = {521–557}
}

@article{Evgrafov2014,
title = {State space {Newton’s} method for topology optimization},
journal = {Computer Methods in Applied Mechanics and Engineering},
volume = {278},
pages = {272-290},
year = {2014},
issn = {0045-7825},
doi = {https://doi.org/10.1016/j.cma.2014.06.005},
url = {https://www.sciencedirect.com/science/article/pii/S004578251400190X},
author = {Evgrafov, A.}
}

@book{allgower2012numerical,
  title={Numerical continuation methods: an introduction},
  author={Allgower, E. L. and Georg, K.},
  volume={13},
  year={2012},
  publisher={Springer Science \& Business Media}
}

@article{Malinen2010,
	author = {Malinen, I. and Tanskanen, J.},
	doi = {10.1016/j.compchemeng.2010.03.013},
	issn = {0098-1354},
	journal = {Computers \& Chemical Engineering},
	number = {11},
	pages = {1761–1774},
	title = {Homotopy parameter bounding in increasing the robustness of homotopy continuation methods in multiplicity studies},
	url = {https://www.sciencedirect.com/science/article/pii/S0098135410001043},
	volume = {34},
	year = {2010}
}

@article{Dunlavy2005,
title = {Homotopy optimization methods for global optimization.},
author = {Dunlavy, D. M. and O'Leary, D. P.},
abstractNote = {We define a new method for global optimization, the Homotopy Optimization Method (HOM). This method differs from previous homotopy and continuation methods in that its aim is to find a minimizer for each of a set of values of the homotopy parameter, rather than to follow a path of minimizers. We define a second method, called HOPE, by allowing HOM to follow an ensemble of points obtained by perturbation of previous ones. We relate this new method to standard methods such as simulated annealing and show under what circumstances it is superior. We present results of extensive numerical experiments demonstrating performance of HOM and HOPE.},
doi = {10.2172/876373},
url = {https://www.osti.gov/biblio/876373}, journal = {},
number = {} ,
volume = {},
place = {United States},
year = {2005},
month = {12}
}

@misc{CesaranoEndtmayerGangl2024PRE,
      title={Homotopy methods for higher order shape optimization: A globalized shape-{Newton} method and Pareto-front tracing},
      author={A. Cesarano and B. Endtmayer and P. Gangl},
      year={2024},
      eprint={2405.03421},
      archivePrefix={arXiv},
      primaryClass={math.NA},
      url={https://arxiv.org/abs/2405.03421}
}

@article{Forsgren2002IMN,
  author =       "Forsgren, A. and Gill, P. E. and Wright, M. H.",
  title =        "Interior Methods for Nonlinear Optimization",
  journal =      {SIAM review},
  volume =       "44",
  number =       "4",
  pages =        "525--597",
  year =         "2002",
  CODEN =        "SIREAD",
  DOI =          "https://doi.org/10.1137/S0036144502414942",
  ISSN =         "0036-1445 (print), 1095-7200 (electronic)",
  ISSN-L =       "0036-1445",
  MRclass =      "90C51 (65K05 90-02 90C30)",
  MRnumber =     "1980444 (2004c:90098)",
  MRreviewer =   "E. Alper Y{\i}ld{\i}r{\i}m",
  bibdate =      "Tue Oct 22 18:27:55 MDT 2002",
  bibsource =    "http://epubs.siam.org/sam-bin/dbq/toc/SIREV/44/4;
                 http://www.math.utah.edu/pub/bibnet/authors/f/forsgren-anders.bib",
  URL =          "http://epubs.siam.org/sam-bin/dbq/article/41494",
  acknowledgement = ack-nhfb,
  fjournal =     "SIAM Review",
  journal-URL =  "http://epubs.siam.org/sirev"
}

@article{Villalobos2004,
  title = {Sphere of Convergence of {Newton’s} Method on Two Equivalent Systems from Nonlinear Programming},
  volume = {121},
  ISSN = {1573-2878},
  url = {http://dx.doi.org/10.1023/B:JOTA.0000037601.54325.3d},
  DOI = {10.1023/b:jota.0000037601.54325.3d},
  number = {3},
  journal = {Journal of Optimization Theory and Applications},
  publisher = {Springer Science and Business Media LLC},
  author = {Villalobos,  M. C. and Tapia,  R. A. and Zhang,  Y.},
  year = {2004},
  month = jun,
  pages = {489-514}
}

@book{Nocedal2006NOpt,
  address = {New York, NY},
  author = {Nocedal, {J.} and Wright, {S. J.}},
  edition = {2. ed.},
  isbn = {978-0-387-30303-1},
  pagetotal = {XXII, 664},
  ppn_gvk = {502988711},
  publisher = {Springer},
  series = {Springer series in operations research and financial engineering},
  title = {Numerical optimization},
  url = {http://gso.gbv.de/DB=2.1/CMD?ACT=SRCHA&SRT=YOP&IKT=1016&TRM=ppn+502988711&sourceid=fbw_bibsonomy},
  year = 2006
}

@article{Sigmund1998,
  title = {Numerical instabilities in topology optimization: A survey on procedures dealing with checkerboards,  mesh-dependencies and local minima},
  volume = {16},
  ISSN = {1615-1488},
  url = {http://dx.doi.org/10.1007/BF01214002},
  DOI = {10.1007/bf01214002},
  number = {1},
  journal = {Structural Optimization},
  publisher = {Springer Science and Business Media LLC},
  author = {Sigmund,  O. and Petersson,  J.},
  year = {1998},
  month = aug,
  pages = {68-75}
}

@article{Modica1987,
  title = {The gradient theory of phase transitions and the minimal interface criterion},
  volume = {98},
  ISSN = {1432-0673},
  url = {http://dx.doi.org/10.1007/BF00251230},
  DOI = {10.1007/bf00251230},
  number = {2},
  journal = {Archive for Rational Mechanics and Analysis},
  publisher = {Springer Science and Business Media LLC},
  author = {Modica,  L.},
  year = {1987},
  month = jun,
  pages = {123–142}
}

@article{schmidt2008pareto,
author = {Schmidt, S. and Schulz, V.},
year = {2008},
month = {05},
pages = {},
title = {Pareto-curve continuation in multi-objective optimization},
volume = {4},
journal = {Pacific Journal of Optimization}
}

@article{MartinSchuetze2018,
  title = {Pareto Tracer: a predictor–corrector method for multi-objective optimization problems},
  volume = {50},
  ISSN = {1029-0273},
  url = {http://dx.doi.org/10.1080/0305215X.2017.1327579},
  DOI = {10.1080/0305215x.2017.1327579},
  number = {3},
  journal = {Engineering Optimization},
  publisher = {Informa UK Limited},
  author = {Martín,  A. and Sch\"{u}tze,  O.},
  year = {2017},
  month = jun,
  pages = {516–536}
}

@book{Schwetlick1979,
  title={Numerische L{\"o}sung nichtlinearer Gleichungen},
  author={Schwetlick, H.},
  series={Mathematik f{\"u}r Naturwissenschaft und Technik},
  url={https://books.google.at/books?id=TcYPzQEACAAJ},
  year={1979},
  publisher={VEB Deutscher Verlag der Wisssenschaften}
}
\bibliographystyle{plain}
\end{document}